\RequirePackage{fix-cm}
\documentclass[smallextended]{svjour3}       % onecolumn (second format)
\smartqed  % flush right qed marks, e.g. at end of proof
\usepackage{graphicx}
\usepackage{amssymb}
\usepackage{amsmath}
\usepackage{tikz}
\usepackage{textcomp}
\journalname{Complex Analysis and Operator Theory}
\begin{document}

\title{Application of the Weyl calculus perspective on discrete octonionic analysis in bounded domains%\thanks{Grants or other notes
}

\titlerunning{Weyl approach to discrete octonionic analysis in bounded domains}        % if too long for running head

\author{Rolf S\"oren Krau\ss har \and Anastasiia Legatiuk \and Dmitrii Legatiuk}

\authorrunning{R.S. Krau\ss har, A. Legatiuk, D. Legatiuk} % if too long for running head

\institute{Rolf S\"oren Krau\ss har \at
	       University of Erfurt \\
		   Chair of Mathematics \\
           Nordh\"auser Str. 63, 99089 Erfurt, Germany \\
           \email{soeren.krausshar@uni-erfurt.de}    
           \and
		   Dmitrii Legatiuk \at
           University of Erfurt \\
		   Chair of Mathematics \\
           Nordh\"auser Str. 63, 99089 Erfurt, Germany \\
           \email{dmitrii.legatiuk@uni-erfurt.de}
		   \and
		   Anastasiia Legatiuk \at
           University of Erfurt \\
		   Chair of Mathematics \\
           Nordh\"auser Str. 63, 99089 Erfurt, Germany \\
           \email{anastasiia.legatiuk@gmail.com}
}

\date{Received: date / Accepted: date}
\dedication{This paper is dedicated to the memory of Frank Sommen.}

\maketitle

\begin{abstract}

In this paper, we finish the basic development of the discrete octonionic analysis by presenting a Weyl calculus-based approach to bounded domains in $\mathbb{R}^{8}$. In particular, we explicitly prove the discrete Stokes formula for a bounded cuboid, and then we generalise this result to arbitrary bounded domains in interior and exterior settings by the help of characteristic functions. After that, discrete interior and exterior Borel-Pompeiu and Cauchy formulae are introduced. Finally, we recall the construction of discrete octonionic Hardy spaces for bounded domains. Moreover, we explicitly explain where the non-associativity of octonionic multiplication is essential and where it is not. Thus, this paper completes the basic framework of the discrete octonionic analysis introduced in previous papers, and, hence, provides a solid foundation for further studies in this field.

%\PACS{PACS code1 \and PACS code2 \and more}
\end{abstract}

\section{Introduction}

The set of octonions represents an eight-dimensional generalisation of complex numbers. Consider $8$-dimensional Euclidean space $\mathbb{R}^{8}$ with the basis unit vectors $\mathbf{e}_{k}$, $k=0,1,\ldots,7$ and points $\mathbf{x}=(x_{0},x_{1},\ldots, x_{7})$. Now an element of $\mathbb{R}^{8}$ can be expressed as follows
\begin{equation*}
x = x_{0}\mathbf{e}_{0}+x_{1}\mathbf{e}_{1}+x_{2}\mathbf{e}_{2}+x_{3}\mathbf{e}_{3}+x_{4}\mathbf{e}_{4}+x_{5}\mathbf{e}_{5}+x_{6}\mathbf{e}_{6}+x_{7}\mathbf{e}_{7},
\end{equation*}
where the basis elements additionally satisfy the following rules  $\mathbf{e}_{4}=\mathbf{e}_{1}\mathbf{e}_{2}$, $\mathbf{e}_{5}=\mathbf{e}_{1}\mathbf{e}_{3}$, $\mathbf{e}_{6}=\mathbf{e}_{2}\mathbf{e}_{3}$ and $\mathbf{e}_{7}=\mathbf{e}_{4}\mathbf{e}_{3}=(\mathbf{e}_{1}\mathbf{e}_{2})\mathbf{e}_{3}$. Moreover, we have $\mathbf{e}_{i}^{2}=-1$ and  $\mathbf{e}_{0}\mathbf{e}_{i}=\mathbf{e}_{i}\mathbf{e}_{0}$ for all $i=1,\ldots,7$, and $\mathbf{e}_{i}\mathbf{e}_{j}=-\mathbf{e}_{j}\mathbf{e}_{i}$ for all mutually distinct $i,j\in\left\{1,\ldots,7\right\}$. The particular feature of octonions is that their multiplication operation is closed but not associative, because we have $(\mathbf{e}_{i}\mathbf{e}_{j})\mathbf{e}_{k}=-\mathbf{e}_{i}(\mathbf{e}_{j}\mathbf{e}_{k})$ for mutually distinct and non-zero $i,j,k$ and if furthermore $\mathbf{e}_{i}\mathbf{e}_{j}\neq \pm \mathbf{e}_k$. Because of the non-associativity, octonionic analysis has not particularly been well studied for a long time. However in the recent years, there is a growing interest in the development of consistent function theories in the octonionic setting, see for example \cite{ConKra2021,FL,Kauhanen_1,QR2022,WL2018} as well as the references therein.\par
Looking at practical applications of function theories and Clifford analysis, it is clear that the continuous constructions need to be discretised for supporting more efficient numerical solutions. The drawback of a direct discretisation of the continuous constructions is that all important quantities of PDEs and operators are only approximated and not represented exactly at the discrete level. To overcome this drawback, construction of discrete counterparts of classical (associative) continuous Clifford analysis has been an area of active research for many years. In particular, Frank Sommen has contributed significantly to this field by his groundbreaking paper \cite{Sommen_1}.\par
Originated by the pioneering work of J. Vaz in 1997 \cite{Vaz}, the discrete Clifford analysis has been intensively developed by Frank Sommen and his co-authors in 2009-2010 in a series of papers \cite{Brackx_1,Frank_1,Frank_2}. In these works, the basic framework of, what is nowadays referred to as discrete Clifford analysis, rendering some even earlier ideas from \cite{FKS07,Guerlebeck_3} has been presented. The main advantage of the discretisation based on the Weyl relation is that we preserve the classical factorisation of the star-Laplacian by the help of two discrete Dirac or Cauchy-Riemann operators. This factorisation opens many possibilities to use this discrete setting in various situations. In particular, Frank Sommen and his co-authors have also addressed discrete heat equations \cite{Baaske}, Taylor expansions \cite{De_Ridder_1}, rotations \cite{De_Ridder_2}, and discrete Hardy spaces \cite{CKKS}. Particularly the definition of the discrete Hardy spaces allowed then to study boundary values of discrete monogenic functions. These results have been extended to bounded domains in \cite{Cerejeiras_2,Cerejeiras}. It is also worth mentioning the old paper of Frank Sommen from 1984 on monogenic differential forms and homology theory \cite{Frank_3}, which, although being not directly related to the discrete setting, inspired thirty years later a completely new direction of research within the Clifford algebra community -- script geometry \cite{Cerejeiras_3}, which can also lead to numerical schemes, however not based on finite differences, but on finite elements \cite{Cerejeiras_4}. It is worth mentioning that the script geometry might be seen as a specific version of Whitney stratified chains and cochains \cite{Goresky}, but the essential difference between the script geometry and this theory or similar theories is the concept of tightness. This concept is purely graph-theoretical and replaces the necessity to study Poincaré lemma and related fundamental constructions, which is not the case in other theories.\par
As we see from the discussion above, the associative setting of discrete Clifford analysis has been an area of active research strongly influenced by Frank Sommen and his ideas. However, the non-associative octonionic case has not been addressed in the discrete setting until the recent time. First ideas on the development of a discrete counterpart of octonionic analysis have been presented in \cite{Krausshar_1,Krausshar_2}, where the continuous octonionic Cauchy-Riemann operators have been discretised by the help of forward and backward finite difference operators. Therefore, an octonionic {\itshape discrete forward Cauchy-Riemann operator} $D_{h}^{+}$ and an octonionic {\itshape discrete backward Cauchy-Riemann operator} $D_{h}^{-}$ have been introduced
\begin{equation*}
D^{+}_{h}:=\sum_{j=0}^{7} \mathbf{e}_j\partial_{h}^{+j}, \quad D^{-}_{h}:=\sum_{j=0}^{7} \mathbf{e}_j\partial_{h}^{-j}.
\end{equation*}
In \cite{Krausshar_2}, a discrete function theory including discrete Borel-Pompeiu formula and discrete Hardy spaces has been introduced for these operators. In \cite{Ren}, ideas of working with bounded domains within the setting from \cite{Krausshar_1,Krausshar_2} have been presented. The drawback of working with the operators $D_{h}^{+}$ and $D_{h}^{-}$ is that they do not factorise the classical star-Laplacian $\Delta_{h}:=\sum\limits_{j=0}^{7}\partial_{h}^{+j}\partial_h^{-j}$ in the canonical form, but rather lead to the more complicated factorisation:
\begin{equation*}
\Delta_{h} = \frac{1}{2}\left(D_{h}^{+}\overline{D_{h}^{-}}+D_{h}^{-}\overline{D_{h}^{+}}\right), 
\end{equation*}
where $\overline{D_{h}^{-}}$ and $\overline{D_{h}^{+}}$ are conjugated operators. This factorisation problem is well-known also in the associative case of Clifford analysis, see e.g. \cite{Brackx_1}. For keeping the canonical factorisation of the star-Laplacian it is necessary to work with the ideas inspired by Frank Sommen and switch to the setting of discrete Clifford analysis.\par
First ideas to use the Weyl calculus approach to discretise the non-associative octonionic setting have been presented in \cite{Krausshar_3}, where a discrete Stokes' formula has been proved. Moreover, it has been shown, that in addition to the classical non-commutative relation used in the Weyl setting, extra non-associative conditions must be posed on the splitting of the basis elements $\mathbf{e}_{i}$ into ``positive'' and ``negative'' parts. However, the fact that the octonions contain also an associative part has been overseen during this construction, as it was pointed out in \cite{Krausshar_4}. Therefore, this paper serves two purposes: at first, we revisit and further develop the Weyl calculus-based approach to discrete octonionic analysis for half-spaces, and after that, we complement the construction by studying bounded domains in $\mathbb{R}^{8}$. The case of bounded domains within the Weyl calculus setting is the
only missing gap in the set up of the general basics of discrete octonionic function theory systematically developed in \cite{Krausshar_1,Krausshar_2,Krausshar_3,Krausshar_4,Ren}. Thus, this paper finishes the basic development of the discrete octonionic function theory, and, hence, along with other works it provides a sound foundation for further research on deeper topics in the discrete octonionic analysis setting.\par

\section{Preliminaries and notations}\label{Section:Preliminaries}

\subsection{Geometrical setting}

Let us consider $8$-dimensional Euclidean space $\mathbb{R}^{8}$ with the standard basis unit vectors $\mathbf{e}_{k}$, $k=0,1,\ldots,7$ and points $\mathbf{x}=(x_{0},x_{1},\ldots, x_{7})$. For $h>0$, we establish over $\mathbb{R}^{8}$ the unbounded uniform lattice $h\mathbb{Z}^{8}$:
\begin{equation*}
h \mathbb{Z}^{8} :=\left\{\mathbf{x} \in {\mathbb{R}}^{8}\,|\, \mathbf{x} = (m_{0}h, m_{1}h,\ldots, m_{7}h), m_{j} \in \mathbb{Z}, j=0,1,\ldots,7\right\}.
\end{equation*}
Additionally, we will also consider upper and lower half-spaces (or half-lattices), which are defined as follows
\begin{equation*}
\begin{array}{rcl}
h\mathbb{Z}_{+}^{8} & := & \left\{(h\underline{m},hm_{7})\colon \underline{m}\in\mathbb{Z}^{7},m_{7}\in\mathbb{Z}_{+}\right\}, \\
h\mathbb{Z}_{-}^{8} & := & \left\{(h\underline{m},hm_{7})\colon \underline{m}\in\mathbb{Z}^{7},m_{7}\in\mathbb{Z}_{-}\right\}.
\end{array}
\end{equation*}\par
Let now $\Omega\subset\mathbb{R}^{8}$ be a bounded, simply connected domain with piecewise smooth boundary $\partial\Omega$. The discrete bounded domain $\Omega_{h}$ is then constructed as follows
\begin{equation*}
\Omega_{h} := \Omega\cap h\mathbb{Z}^{8} = \left\{ mh = \left(m_{0}h,m_{1}h,\ldots,m_{7}h\right) \mid mh\in \Omega\cap h\mathbb{Z}^{8}\right\}
\end{equation*}\par
While dealing with discrete function theories in bounded domains, it is crucially important to provide a clear characterisation of a discrete boundary. In particular, unlike in the continuous case, the bounded setting in the discrete case requires the consideration of a three-layer boundary, where the middle layer can be seen as a \textquotedblleft true\textquotedblright\, boundary, and two other layers are required for studying interior and exterior problems. Several approaches to the description of the discrete boundary have been presented in the past. A constructive approach based on the explicit construction of boundary layers has been used in \cite{Guerlebeck_3,Guerlebeck_10,Legatiuk_1}. The advantage of this approach is that an explicit discretisation algorithm for practical calculations can be easily derived, see \cite{Guerlebeck_8}. Evidently, this approach becomes more bulky with growing dimension. For example in \cite{Cerejeiras_2}, the constructive approach has been used to develop discrete Hardy spaces for bounded domains in $\mathbb{R}^{3}$, which has led to more complicated expressions and constructions. To overcome the growing complexity of working with the constructive approach, a more abstract perspective on the description of the discrete boundary has been presented in \cite{Cerejeiras}. In this case, a topological perspective on sets is adapted to the discrete setting, and the boundary layers are characterised by the help of characteristic functions. The abstract approach to the discrete boundary is arguable more elegantly from the theoretical perspective. However, if we are interested in solving boundary value problems in practice, then we need an interplay between both constructive and abstract approaches. Moreover, for the purpose of developing discrete octonionic analysis for bounded domains, it is necessary to \textquotedblleft catch\textquotedblright\, the effect of non-associativity of octonionic multiplication, which might be hidden by the abstract approach to the discretisation. Therefore, in this paper, we will mix both approaches to the discretisation and switch between them by convenience reasons.\par
Following \cite{Cerejeiras}, we now introduce the following topological characterisation of discrete domains:
\begin{definition}
For a given discrete domain $\Omega_{h}$, the following objects are introduced:
\begin{itemize}
\item[(i)] a discrete complementary domain to $\Omega_{h}$: $\Omega_{h}^c := h \mathbb{Z}^{8}\setminus \Omega_{h};$
\item[(ii)] the discrete interior of $\Omega_{h}$, denoted by $int(\Omega_{h})$, is the set of all points $mh\in \Omega_{h}$ such that at least one of its immediate neighbour points $(m + \mathbf{e}_{j})h, (m - \mathbf{e}_{j})h$ for some $j=0,\ldots,7$, also belongs to $\Omega_{h}$, i.e.
\begin{equation*}
int(\Omega_{h}) := \left\{ mh \in \Omega_{h}\,|\, \exists\, j \colon (m + \mathbf{e}_{j})h \in \Omega_{h} \vee (m - \mathbf{e}_{j})h \in \Omega_{h}\right\};
\end{equation*}
\item[(iii)] the discrete exterior of $\Omega_{h}$, denoted by $\Omega_h^{ext}$, is defined symmetrically as the interior of the complementary domain $\Omega_{h}^c$.
\end{itemize}
\end{definition}\par
Now we can provide a classification of boundary points into three boundary layers: 
\begin{definition}\label{Definition_geometry}
Let $\Omega_{h}$ be a discrete domain $h\mathbb{Z}^{8}$ with the lattice constant $h>0.$ We say that
\begin{itemize}
\item[(i)] a point $mh \in \Omega_{h}$ is a point of the {\itshape interior boundary layer} $\gamma_{h}^{+}$, if at least one of its neighbour points does not belong to $\Omega_{h}$, i.e.
\begin{equation*}
\gamma_{h}^{+} := \left\{mh \in \Omega_{h}\,|\, \exists \,j \colon (m + \mathbf{e}_{j})h \notin \Omega_{h} \vee (m - \mathbf{e}_{j})h \notin \Omega_{h} \right\};
\end{equation*}
\item[(ii)] a point $mh$ is a point of the {\itshape middle boundary layer} $\gamma_{h}^{\ast}$, if its neighbourhood contains points belonging to $\Omega_{h}$, as well as points belonging to $\Omega_h^{ext}$, i.e.
\begin{equation*}
\gamma_{h}^{\ast} := \left\{mh \,|\, \exists \,j \colon (m \pm \mathbf{e}_{j})h \in \Omega_{h} \wedge (m \mp \mathbf{e}_{j})h \in \Omega_h^{ext} \right\};
\end{equation*}
\item[(iii)] a point $mh \in \Omega_h^{ext}$ is a point of the {\itshape exterior boundary layer} $\gamma_{h}^{-}$, if at least one of its neighbour points does not belong to $\Omega_h^{ext}$, i.e.
\begin{equation*}
\gamma_{h}^{-} := \left\{mh \in \Omega_h^{ext}\,|\, \exists \,j \colon (m + \mathbf{e}_{j})h \notin \Omega_h^{ext} \vee (m - \mathbf{e}_{j})h \notin \Omega_h^{ext} \right\}.
\end{equation*}
\end{itemize}
\end{definition}\par
It is worth underlining that the middle layer $\gamma_{h}^{\ast}$ can be regarded as a \textquotedblleft true\textquotedblright\, boundary of a domain. Of course, the situation is trickier in the discrete case, because if a domain with curved boundaries is approximated by the lattice $h\mathbb{Z}^{8}$, then it could happen that none of the elements of $\gamma_{h}^{\ast}$ indeed belongs to the continuous boundary $\partial\Omega$. But this effect will be eliminated by considering the convergence process for $h\to 0$.\par
Next, we introduce the characteristic functions for the discrete domains $\Omega_{h}$ and $\Omega_{h}^{ext}$ as follows
\begin{equation*}
\chi_{\Omega_{h}}(mh) := \left\{\begin{array}{cl}
1, &  mh \in \Omega_{h}, \\
0, & \mbox{otherwise},
\end{array} \right. \qquad \chi_{\Omega_{h}^{ext}} (mh) := \left\{\begin{array}{cl}
1, & mh \in \Omega_{h}^{ext}, \\
0, & \mbox{otherwise}.
\end{array} \right.
\end{equation*}
The main advantage of working with characteristic functions is their property that the forward and backward derivatives (introduced right below) applied to them vanish everywhere except the boundary layers. This property admits another characterisation of boundary layers, as well as a shorter presentation of otherwise bulky expressions, see \cite{Cerejeiras} for details.\par 
Let us now introduce the classical forward and backward differences $\partial_{h}^{\pm j}$
\begin{equation}
\label{Finite_differences}
\begin{array}{lcl}
\partial_{h}^{+j}f(mh) & := & h^{-1}(f(mh+\mathbf{e}_jh)-f(mh)), \\
\partial_{h}^{-j}f(mh) & := & h^{-1}(f(mh)-f(mh-\mathbf{e}_jh)),
\end{array}
\end{equation}
for discrete functions $f(mh)$ with $mh\in h\mathbb{Z}^{8}$. Moreover, as usual, we consider discrete functions $f\colon \Omega_{h} \subset  h\mathbb{Z}^{8} \to \mathbb{O}$ in this paper, and all important properties of such functions, for example the $l^{p}$-summability ($1\leq p<\infty$), are defined component-wisely.\par
\begin{remark}
We would like to remark, that the discrete setting has the advantage, that there is no need to introduce sophisticated expressions for normal vectors, because normal vectors for the lattice $h\mathbb{Z}^{8}$ and discrete domains are always coincide with the basis elements $\mathbf{e}_{j}$, $j=0,1,\ldots,7$ and their negative counterparts. 
\end{remark}\par

\subsection{Weyl calculus approach}

In this short section, we briefly recall basic ideas of the Weyl calculus approach to the discretisation of octonionic analysis. Therefore, following the classical ideas from \cite{Brackx_1,FKS07}, each basis element $\mathbf{e}_{k}$, $k=0,1,\ldots,7$, is sub-divided into positive and negative directions $\mathbf{e}_{k}^+$ and $\mathbf{e}_{k}^-$, $k=0,1,\ldots,7$, i.e., $\mathbf{e}_{k}=\mathbf{e}_{k}^{+}+\mathbf{e}_{k}^-$. Moreover, this splitting must satisfy the following relations
\begin{equation}
\label{Splitting_relations}
\left\{
\begin{array}{ccc}
\mathbf{e}_j^{-}\mathbf{e}_k^{-}+\mathbf{e}_k^{-}\mathbf{e}_j^{-} &=&0, \\
\mathbf{e}_j^{+}\mathbf{e}_k^{+}+\mathbf{e}_k^{+}\mathbf{e}_j^{+}&=&0, \\
\mathbf{e}_j^{+}\mathbf{e}_k^{-} + \mathbf{e}_k^{-}\mathbf{e}_j^{+}&=&-\delta_{jk},
\end{array}
\right.
\end{equation}
where $\delta_{jk}$ is the Kronecker delta symbol. Moreover, taking into account that octonionic multiplication is not associative,
\begin{equation*}
(\mathbf{e}_{i}\mathbf{e}_{j})\mathbf{e}_{k}=-\mathbf{e}_{i}(\mathbf{e}_{j}\mathbf{e}_{k}),
\end{equation*}
it is necessary to establish extra conditions for the splitting of the basis elements. In \cite{Krausshar_3} it has been shown, that the following anti-associative relations must be satisfied: 
\begin{equation}
\label{Splitting_anti_associative_relations}
\begin{array}{rclrcl}
\left(\mathbf{e}_{i}^{+}\mathbf{e}_{j}^{+}\right)\mathbf{e}_{k}^{+} & = & -\mathbf{e}_{i}^{+}\left(\mathbf{e}_{j}^{+}\mathbf{e}_{k}^{+}\right), & \left(\mathbf{e}_{i}^{+}\mathbf{e}_{j}^{+}\right)\mathbf{e}_{k}^{-} & = & -\mathbf{e}_{i}^{+}\left(\mathbf{e}_{j}^{+}\mathbf{e}_{k}^{-}\right), \\
\left(\mathbf{e}_{i}^{-}\mathbf{e}_{j}^{+}\right)\mathbf{e}_{k}^{+} & = & -\mathbf{e}_{i}^{-}\left(\mathbf{e}_{j}^{+}\mathbf{e}_{k}^{+}\right), & \left(\mathbf{e}_{i}^{-}\mathbf{e}_{j}^{+}\right)\mathbf{e}_{k}^{-} & = & -\mathbf{e}_{i}^{-}\left(\mathbf{e}_{j}^{+}\mathbf{e}_{k}^{-}\right), \\
\left(\mathbf{e}_{i}^{+}\mathbf{e}_{j}^{-}\right)\mathbf{e}_{k}^{+} & = & -\mathbf{e}_{i}^{+}\left(\mathbf{e}_{j}^{-}\mathbf{e}_{k}^{+}\right), & \left(\mathbf{e}_{i}^{+}\mathbf{e}_{j}^{-}\right)\mathbf{e}_{k}^{-} & = & -\mathbf{e}_{i}^{+}\left(\mathbf{e}_{j}^{-}\mathbf{e}_{k}^{-}\right), \\
\left(\mathbf{e}_{i}^{-}\mathbf{e}_{j}^{-}\right)\mathbf{e}_{k}^{+} & = & -\mathbf{e}_{i}^{-}\left(\mathbf{e}_{j}^{-}\mathbf{e}_{k}^{+}\right), & \left(\mathbf{e}_{i}^{-}\mathbf{e}_{j}^{-}\right)\mathbf{e}_{k}^{-} & = & -\mathbf{e}_{i}^{-}\left(\mathbf{e}_{j}^{-}\mathbf{e}_{k}^{-}\right).
\end{array}
\end{equation}
It is necessary to remark, that these non-associative conditions are satisfied only for mutually distinct and non-zero indices $i,j,k$ satisfying additionally $\mathbf{e}_{i}\mathbf{e}_{j}\neq \pm \mathbf{e}_k$, otherwise, the multiplication is associative $(\mathbf{e}_{i}\mathbf{e}_{j})\mathbf{e}_{k}=\mathbf{e}_{i}(\mathbf{e}_{j}\mathbf{e}_{k})$. This aspect of the octonionic multiplication has been overseen in \cite{Krausshar_3}. Therefore, in Section~\ref{Section:Half_spaces}, we will revisit this setting and present corrected formulae, before addressing the case of bounded domains.\par

\subsection{Discrete octonionic monogenicity and discrete fundamental solutions}

In this short section, we briefly recall ideas on discrete octonionic monogenicity in the Weyl-calculus setting and some basic facts about the discrete fundamental solutions of the discrete Cauchy-Riemann operators, see \cite{Krausshar_4} for more details. Let $\Omega_{h}\subset h\mathbb{Z}^{8}$, then the discrete Cauchy-Riemann operator $D^{+-}_{h}\colon l^{p}(\Omega_{h},\mathbb{O})\to l^{p}(\Omega_{h},\mathbb{O})$ and its adjoint operator $D^{-+}_{h}\colon l^{p}(\Omega_{h},\mathbb{O})\to l^{p}(\Omega_{h},\mathbb{O})$ are introduced as follows
\begin{equation*}
D^{+-}_{h}:=\sum_{j=0}^{7} \mathbf{e}_{j}^{+}\partial_{h}^{+j}+\mathbf{e}_{j}^{-}\partial_{h}^{-j}, \quad D^{-+}_{h}:=\sum_{j=0}^{7} \mathbf{e}_{j}^{+}\partial_{h}^{-j}+\mathbf{e}_{j}^{-}\partial_{h}^{+j}.
\end{equation*}
These operators preserve the canonical factorisation of the star-Laplacian 
\begin{equation*}
(D^{+-}_{h})^2=(D^{-+}_h)^2 = -\Delta_{h},
\end{equation*}
with
\begin{equation*}
\Delta_{h}:=\sum_{j=0}^{7}\partial_{h}^{+j}\partial_h^{-j}.
\end{equation*}\par
By the help of the discrete Cauchy-Riemann operators, {\itshape discrete  monogenic functions} can now be introduced:
\begin{definition}
A function $f\in l^{p}(\Omega_{h},\mathbb{O})$ is called {\itshape discrete left monogenic} if $D_{h}^{+-}f=0$ in $\Omega_{h}$. Respectively, a function $f\in l^{p}(\Omega_{h},\mathbb{O})$ is called {\itshape discrete left anti-monogenic} if $D_{h}^{-+}f=0$ in $\Omega_{h}$.
\end{definition}\par
Finally, we recall the definition of a discrete fundamental solution of the discrete Cauchy-Riemann operator:
\begin{definition}\label{Definition:Discrete_fundamental_solution}
The function $E_{h}^{+-}\colon h\mathbb{Z}^{8} \rightarrow \mathbb{O}$ is called a {\itshape discrete fundamental solution} of $D_{h}^{+-}$ if it satisfies
\begin{equation*}
D_{h}^{+-}E_{h}^{+-} =\delta_h = \begin{cases} h^{-8}, & \mbox{for } mh=0,\\ 
0,& \mbox{for } mh\neq 0,
\end{cases}
\end{equation*} 
for all grid points $mh$ of $h\mathbb{Z}^{8}$. Analogously, The function $E_{h}^{-+}\colon h\mathbb{Z}^{8} \rightarrow \mathbb{O}$ is called a {\itshape discrete fundamental solution} of $D_{h}^{-+}$ if it satisfies
\begin{equation*}
D_{h}^{-+}E_{h}^{-+} =\delta_h = \begin{cases} h^{-8}, & \mbox{for } mh=0,\\ 
0,& \mbox{for } mh\neq 0,
\end{cases}
\end{equation*} 
for all grid points $mh$ of $h\mathbb{Z}^{8}$.
\end{definition}
We will only present the final form of the discrete fundamental solutions $E_{h}^{+-}$ and $E_{h}^{-+}$, and we refer to \cite{CKKS} for the associative setting of discrete Clifford analysis and to \cite{Krausshar_2} for the discrete octonionic setting. At first, let us recall the Fourier symbols of the operators $D_{h}^{+-}$ and $D_{h}^{-+}$:
\begin{equation*}
\tilde{\xi}_{+} = \sum\limits_{j=0}^{7} \left(\mathbf{e}_{j}^{+}\xi_{h}^{+j} + \mathbf{e}_{j}^{-}\xi_{h}^{-j}\right) \mbox{ and } \tilde{\xi}_{-} = \sum\limits_{j=0}^{7} \left(\mathbf{e}_{j}^{+}\xi_{h}^{-j} + \mathbf{e}_{j}^{-}\xi_{h}^{+j}\right),
\end{equation*}
respectively, where $\xi_{h}^{\pm j}=\mp h^{-1}\left( 1-e^{\mp ih\xi_j}\right)$. These symbols are used to express the discrete fundamental solutions as follows
\begin{equation}
\label{Discrete_fundamental_solutions}
\begin{array}{ccl}
\displaystyle E_{h}^{+-} & = & \displaystyle \mathcal{R}_{h}\mathcal{F} \left(\frac{\widetilde{\xi}_{+}}{d^2}\right)
= \sum\limits_{j=0}^{7}\mathbf{e}_{j}^{+}\mathcal{R}_{h} \mathcal{F} \left( \frac{\xi_{h}^{+j}}{d^2} \right) + \mathbf{e}_{j}^{-}\mathcal{R}_{h} \mathcal{F} \left( \frac{\xi_{h}^{-j}}{d^2} \right), \\
\displaystyle E_{h}^{-+} & = & \displaystyle \mathcal{R}_{h}\mathcal{F} \left(\frac{\widetilde{\xi}_{-}}{d^2}\right)
= \sum\limits_{j=0}^{7}\mathbf{e}_{j}^{+}\mathcal{R}_{h} \mathcal{F} \left( \frac{\xi_{h}^{-j}}{d^2} \right) + \mathbf{e}_{j}^{-}\mathcal{R}_{h} \mathcal{F} \left( \frac{\xi_{h}^{+j}}{d^2} \right),
\end{array}
\end{equation}
where $d^{2}=\frac{4}{h^{2}}\sum\limits_{j=0}^{7}\sin^{2}\left(\frac{\xi_{j}h}{2}\right)$ is the symbol of the star-Laplacian, $\mathcal{F}$ is the classical continuous Fourier transform
\begin{equation*}
x \mapsto \mathcal{F}f(x) = \frac{1}{(2\pi)^{8}}\int_{\mathbb{R}^{8}} e^{-i \langle x,\xi \rangle} f(\xi)d\xi,
\end{equation*}
applied to an octonionic-valued function $f \in l^{p}\left( h\mathbb{Z}^8,\mathbb{O}\right)$ with $\mathrm{supp}\, f \in\left[-\frac{\pi}{h},\frac{\pi}{h}\right]^{8}$, and $\mathcal{R}_{h}$ is its restriction to the lattice $h\mathbb{Z}^{8}$.\par

\section{Discrete octonionic function theory for half-spaces based on Weyl calculus approach}\label{Section:Half_spaces}

In this section, we revisit the discrete octonionic function theory for half-spaces by taking into account also the associative sub-part of octonions. As we will see, this consideration will make the associator appear in all of the construction, and, thus, making the results presented in this paper closer to the classical continuous setting.\par
We start with the following theorem presenting a discrete Stokes' formula for the whole lattice:
\begin{theorem}\label{Discrete_Stokes_space}
The discrete Stokes' formula for the whole lattice $h\mathbb{Z}^{8}$ is given by
\begin{equation}
\label{Discrete_Stokes_whole_space}
\begin{array}{c}
\displaystyle \sum_{m\in \mathbb{Z}^{8}}  \left[ \left( g(mh)D_h^{-+}\right) f(mh) + g(mh) \left( D_h^{+-}f(mh) \right) \right] h^8 =  \\
\displaystyle = 2\sum_{m\in \mathbb{Z}^{8}} \sum\limits_{s=1}^{7}  \sum_{i\in I_{s}} \sum_{\stackrel{j\in I_{s}}{j\neq i}}^{7} \sum_{\stackrel{k=1}{k\notin I_{s}}}^{7} \left[g_{i}(mh)\mathbf{e}_{i}\left(\partial_{h}^{+j}\mathbf{e}_{j}^{+}f_{k}(mh)\mathbf{e}_{k}\right) \right. \\
\displaystyle \left. +g_{i}(mh)\mathbf{e}_{i}\left(\partial_{h}^{-j}\mathbf{e}_{j}^{-}f_{k}(mh)\mathbf{e}_{k}\right)\right]h^{8}
\end{array}  
\end{equation}
for all discrete functions $f$ and $g$ such that the series converge, where the index sets $I_{s}$, $s=1,\ldots,7$ are given by
\begin{equation*}
\begin{array}{cclcclcclccl}
I_{1} & := & \left\{1,2,4\right\}, & I_{2} & := & \left\{1,3,5\right\}, & I_{3} & := & \left\{1,6,7\right\}, & I_{4} & := & \left\{2,3,6\right\}, \\
I_{5} & := & \left\{2,5,7\right\}, & I_{6} & := & \left\{3,4,7\right\}, & I_{7} & := & \left\{4,5,6\right\}.
\end{array}
\end{equation*}
\end{theorem}
\begin{proof}
The proof of this theorem in octonionic setting essentially follows the steps presented in \cite{Krausshar_4}. However, for correcting the proof shown in \cite{Krausshar_3}, we will nonetheless present the main steps explicitly. At first the definition of the Cauchy-Riemann operator $D_{h}^{-+}$ is used on the left-hand site.
\begin{equation*}
\begin{array}{c}
\displaystyle \sum\limits_{m\in \mathbb{Z}^{8}} \left[g(mh)D_{h}^{-+}\right]f(mh)h^{8} = \sum\limits_{m\in \mathbb{Z}} \sum\limits_{j=0}^{7} \left[\partial_{h}^{-j}g(mh)\mathbf{e}_{j}^{+}+\partial_{h}^{+j}g(mh)\mathbf{e}_{j}^{-}\right] f(mh) h^{8} \\
= \displaystyle \sum\limits_{m\in \mathbb{Z}^{8}} \sum\limits_{j=0}^{7}\sum \limits_{i=0}^{7} \left[\partial_{h}^{-j}g_{i}(mh)\left(\mathbf{e}_{i}^{+}+\mathbf{e}_{i}^{-}\right)\mathbf{e}_{j}^{+}+\partial_{h}^{+j}g_{i}(mh)\left(\mathbf{e}_{i}^{+}+\mathbf{e}_{i}^{-}\right)\mathbf{e}_{j}^{-}\right] f(mh) h^{8}.
\end{array}
\end{equation*}
Next, we split also the unit vectors of function $f$. After multiplying the result, the following expression containing 8 summands is obtained
\begin{equation*}
\begin{array}{c}
\displaystyle \sum\limits_{m\in \mathbb{Z}^{8}} \sum\limits_{j=0}^{7}\sum \limits_{i=0}^{7} \sum \limits_{k=0}^{7} \left[\partial_{h}^{-j}g_{i}(mh)f_{k}(mh)\left(\mathbf{e}_{i}^{+}\mathbf{e}_{j}^{+}\right)\mathbf{e}_{k}^{+} + \partial_{h}^{-j}g_{i}(mh)f_{k}(mh)\left(\mathbf{e}_{i}^{+}\mathbf{e}_{j}^{+}\right)\mathbf{e}_{k}^{-} \right. \\
\displaystyle + \partial_{h}^{-j}g_{i}(mh)f_{k}(mh)\left(\mathbf{e}_{i}^{-}\mathbf{e}_{j}^{+}\right)\mathbf{e}_{k}^{+} + \partial_{h}^{-j}g_{i}(mh)f_{k}(mh)\left(\mathbf{e}_{i}^{-}\mathbf{e}_{j}^{+}\right)\mathbf{e}_{k}^{-} \\
\displaystyle + \partial_{h}^{+j}g_{i}(mh)f_{k}(mh)\left(\mathbf{e}_{i}^{+}\mathbf{e}_{j}^{-}\right)\mathbf{e}_{k}^{+} + \partial_{h}^{+j}g_{i}(mh)f_{k}(mh)\left(\mathbf{e}_{i}^{+}\mathbf{e}_{j}^{-}\right)\mathbf{e}_{k}^{-} \\
\displaystyle \left. + \partial_{h}^{+j}g_{i}(mh)f_{k}(mh)\left(\mathbf{e}_{i}^{-}\mathbf{e}_{j}^{-}\right)\mathbf{e}_{k}^{+} + \partial_{h}^{+j}g_{i}(mh)f_{k}(mh)\left(\mathbf{e}_{i}^{-}\mathbf{e}_{j}^{-}\right)\mathbf{e}_{k}^{-} \right] h^{8}.
\end{array}
\end{equation*}
At this point we need to address the multiplication of all splitted unit vectors of the form $\left(\mathbf{e}_{i}^{\pm}\mathbf{e}_{j}^{\pm}\right)\mathbf{e}_{k}^{\pm}$ for all possible combination of $+$ and $-$. The goal is to preserve the octonionic multiplication, which is non-associative for mutually distinct, non-zero $i,j,k$ and $\mathbf{e}_{i}\mathbf{e}_{j} \neq \pm \mathbf{e}_{k}$, and associative otherwise. Hence, to keep this construction also valid in the Weyl calculus setting, it is necessary to impose the following anti-associative relations:
\begin{equation*}
\begin{array}{rclrcl}
\left(\mathbf{e}_{i}^{+}\mathbf{e}_{j}^{+}\right)\mathbf{e}_{k}^{+} & = & -\mathbf{e}_{i}^{+}\left(\mathbf{e}_{j}^{+}\mathbf{e}_{k}^{+}\right), & \left(\mathbf{e}_{i}^{+}\mathbf{e}_{j}^{+}\right)\mathbf{e}_{k}^{-} & = & -\mathbf{e}_{i}^{+}\left(\mathbf{e}_{j}^{+}\mathbf{e}_{k}^{-}\right), \\
\left(\mathbf{e}_{i}^{-}\mathbf{e}_{j}^{+}\right)\mathbf{e}_{k}^{+} & = & -\mathbf{e}_{i}^{-}\left(\mathbf{e}_{j}^{+}\mathbf{e}_{k}^{+}\right), & \left(\mathbf{e}_{i}^{-}\mathbf{e}_{j}^{+}\right)\mathbf{e}_{k}^{-} & = & -\mathbf{e}_{i}^{-}\left(\mathbf{e}_{j}^{+}\mathbf{e}_{k}^{-}\right), \\
\left(\mathbf{e}_{i}^{+}\mathbf{e}_{j}^{-}\right)\mathbf{e}_{k}^{+} & = & -\mathbf{e}_{i}^{+}\left(\mathbf{e}_{j}^{-}\mathbf{e}_{k}^{+}\right), & \left(\mathbf{e}_{i}^{+}\mathbf{e}_{j}^{-}\right)\mathbf{e}_{k}^{-} & = & -\mathbf{e}_{i}^{+}\left(\mathbf{e}_{j}^{-}\mathbf{e}_{k}^{-}\right), \\
\left(\mathbf{e}_{i}^{-}\mathbf{e}_{j}^{-}\right)\mathbf{e}_{k}^{+} & = & -\mathbf{e}_{i}^{-}\left(\mathbf{e}_{j}^{-}\mathbf{e}_{k}^{+}\right), & \left(\mathbf{e}_{i}^{-}\mathbf{e}_{j}^{-}\right)\mathbf{e}_{k}^{-} & = & -\mathbf{e}_{i}^{-}\left(\mathbf{e}_{j}^{-}\mathbf{e}_{k}^{-}\right),
\end{array}
\end{equation*}
for mutually distinct, non-zero $i,j,k$ and satisfying $\mathbf{e}_{i}\mathbf{e}_{j} \neq \pm \mathbf{e}_{k}$, while the associativity is imposed for all other combinations of $i,j,k$. This associative part was overseen in \cite{Krausshar_3}. After using these multiplication rules the resulting expression contains associative and non-associative parts. In order to write the associative part as a complete sum over all possible indices $i,j,k$, we add and subtract the non-associative part. This leads to the following expression, which is written again by the help of a triple summation over indices $i,j,k$ as follows
\begin{equation*}
\begin{array}{c}
\displaystyle \sum\limits_{m\in \mathbb{Z}^{8}} \sum\limits_{j=0}^{7}\sum \limits_{i=0}^{7} \sum \limits_{k=0}^{7} \left[\partial_{h}^{-j}g_{i}(mh)f_{k}(mh)\mathbf{e}_{i}^{+}\left(\mathbf{e}_{j}^{+}\mathbf{e}_{k}^{+}\right) + \partial_{h}^{-j}g_{i}(mh)f_{k}(mh)\mathbf{e}_{i}^{+}\left(\mathbf{e}_{j}^{+}\mathbf{e}_{k}^{-}\right) \right. \\
\displaystyle + \partial_{h}^{-j}g_{i}(mh)f_{k}(mh)\mathbf{e}_{i}^{-}\left(\mathbf{e}_{j}^{+}\mathbf{e}_{k}^{+}\right) + \partial_{h}^{-j}g_{i}(mh)f_{k}(mh)\mathbf{e}_{i}^{-}\left(\mathbf{e}_{j}^{+}\mathbf{e}_{k}^{-}\right) \\
\displaystyle + \partial_{h}^{+j}g_{i}(mh)f_{k}(mh)\mathbf{e}_{i}^{+}\left(\mathbf{e}_{j}^{-}\mathbf{e}_{k}^{+}\right) + \partial_{h}^{+j}g_{i}(mh)f_{k}(mh)\mathbf{e}_{i}^{+}\left(\mathbf{e}_{j}^{-}\mathbf{e}_{k}^{-}\right) \\
\displaystyle \left. + \partial_{h}^{+j}g_{i}(mh)f_{k}(mh)\mathbf{e}_{i}^{-}\left(\mathbf{e}_{j}^{-}\mathbf{e}_{k}^{+}\right) + \partial_{h}^{+j}g_{i}(mh)f_{k}(mh)\mathbf{e}_{i}^{-}\left(\mathbf{e}_{j}^{-}\mathbf{e}_{k}^{-}\right) \right] h^{8}.
\end{array}
\end{equation*}
For writing the non-associative part in a shorter form, we introduce the following index sets:
\begin{equation*}
\begin{array}{cclcclcclccl}
I_{1} & := & \left\{1,2,4\right\}, & I_{2} & := & \left\{1,3,5\right\}, & I_{3} & := & \left\{1,6,7\right\}, & I_{4} & := & \left\{2,3,6\right\}, \\
I_{5} & := & \left\{2,5,7\right\}, & I_{6} & := & \left\{3,4,7\right\}, & I_{7} & := & \left\{4,5,6\right\}.
\end{array}
\end{equation*}
By the help of these index sets, we obtain the following expression for the non-associative terms
\begin{equation*}
\begin{array}{c}
\displaystyle
-2\sum_{m\in \mathbb{Z}^{8}} \sum\limits_{s=1}^{7}  \sum_{i\in I_{s}} \sum_{\stackrel{j\in I_{s}}{j\neq i}}^{7} \sum_{\stackrel{k=1}{k\notin I_{s}}}^{7} \left[\partial_{h}^{-j}g_{i}(mh)f_{k}(mh)\mathbf{e}_{i}^{+}\left(\mathbf{e}_{j}^{+}\mathbf{e}_{k}^{+}\right) + \partial_{h}^{-j}g_{i}(mh)f_{k}(mh)\mathbf{e}_{i}^{+}\left(\mathbf{e}_{j}^{+}\mathbf{e}_{k}^{-}\right) \right. \\
\displaystyle + \partial_{h}^{-j}g_{i}(mh)f_{k}(mh)\mathbf{e}_{i}^{-}\left(\mathbf{e}_{j}^{+}\mathbf{e}_{k}^{+}\right) + \partial_{h}^{-j}g_{i}(mh)f_{k}(mh)\mathbf{e}_{i}^{-}\left(\mathbf{e}_{j}^{+}\mathbf{e}_{k}^{-}\right) \\
\displaystyle + \partial_{h}^{+j}g_{i}(mh)f_{k}(mh)\mathbf{e}_{i}^{+}\left(\mathbf{e}_{j}^{-}\mathbf{e}_{k}^{+}\right) + \partial_{h}^{+j}g_{i}(mh)f_{k}(mh)\mathbf{e}_{i}^{+}\left(\mathbf{e}_{j}^{-}\mathbf{e}_{k}^{-}\right) \\
\displaystyle \left. + \partial_{h}^{+j}g_{i}(mh)f_{k}(mh)\mathbf{e}_{i}^{-}\left(\mathbf{e}_{j}^{-}\mathbf{e}_{k}^{+}\right) + \partial_{h}^{+j}g_{i}(mh)f_{k}(mh)\mathbf{e}_{i}^{-}\left(\mathbf{e}_{j}^{-}\mathbf{e}_{k}^{-}\right) \right] h^{8},
\end{array}
\end{equation*}
where the factor 2 comes from the fact, that we added and subtracted the non-associative terms. Next step is to use the definition of finite difference operators~(\ref{Finite_differences}) and perform the change of variables in the resulting expression, see \cite{CKKS,Krausshar_3} for explicit steps. After performing all these steps, we arrive at the expression
\begin{equation*}
\begin{array}{c}
\displaystyle -\sum\limits_{m\in \mathbb{Z}^{8}} \sum\limits_{j=0}^{7}\sum \limits_{i=0}^{7} \sum \limits_{k=0}^{7} \left[\partial_{h}^{+j}g_{i}(mh)f_{k}(mh)\mathbf{e}_{i}^{+}\left(\mathbf{e}_{j}^{+}\mathbf{e}_{k}^{+}\right) + \partial_{h}^{+j}g_{i}(mh)f_{k}(mh)\mathbf{e}_{i}^{+}\left(\mathbf{e}_{j}^{+}\mathbf{e}_{k}^{-}\right) \right. \\
\displaystyle + \partial_{h}^{+j}g_{i}(mh)f_{k}(mh)\mathbf{e}_{i}^{-}\left(\mathbf{e}_{j}^{+}\mathbf{e}_{k}^{+}\right) + \partial_{h}^{+j}g_{i}(mh)f_{k}(mh)\mathbf{e}_{i}^{-}\left(\mathbf{e}_{j}^{+}\mathbf{e}_{k}^{-}\right) \\
\displaystyle + \partial_{h}^{-j}g_{i}(mh)f_{k}(mh)\mathbf{e}_{i}^{+}\left(\mathbf{e}_{j}^{-}\mathbf{e}_{k}^{+}\right) + \partial_{h}^{-j}g_{i}(mh)f_{k}(mh)\mathbf{e}_{i}^{+}\left(\mathbf{e}_{j}^{-}\mathbf{e}_{k}^{-}\right) \\
\displaystyle \left. + \partial_{h}^{-j}g_{i}(mh)f_{k}(mh)\mathbf{e}_{i}^{-}\left(\mathbf{e}_{j}^{-}\mathbf{e}_{k}^{+}\right) + \partial_{h}^{-j}g_{i}(mh)f_{k}(mh)\mathbf{e}_{i}^{-}\left(\mathbf{e}_{j}^{-}\mathbf{e}_{k}^{-}\right) \right] h^{8} +
\end{array}
\end{equation*}
\begin{equation*}
\begin{array}{c}
\displaystyle
+2\sum_{m\in \mathbb{Z}^{8}} \sum\limits_{s=1}^{7}  \sum_{i\in I_{s}} \sum_{\stackrel{j\in I_{s}}{j\neq i}}^{7} \sum_{\stackrel{k=1}{k\notin I_{s}}}^{7} \left[\partial_{h}^{+j}g_{i}(mh)f_{k}(mh)\mathbf{e}_{i}^{+}\left(\mathbf{e}_{j}^{+}\mathbf{e}_{k}^{+}\right) + \partial_{h}^{+j}g_{i}(mh)f_{k}(mh)\mathbf{e}_{i}^{+}\left(\mathbf{e}_{j}^{+}\mathbf{e}_{k}^{-}\right) \right. \\
\displaystyle + \partial_{h}^{+j}g_{i}(mh)f_{k}(mh)\mathbf{e}_{i}^{-}\left(\mathbf{e}_{j}^{+}\mathbf{e}_{k}^{+}\right) + \partial_{h}^{+j}g_{i}(mh)f_{k}(mh)\mathbf{e}_{i}^{-}\left(\mathbf{e}_{j}^{+}\mathbf{e}_{k}^{-}\right) \\
\displaystyle + \partial_{h}^{-j}g_{i}(mh)f_{k}(mh)\mathbf{e}_{i}^{+}\left(\mathbf{e}_{j}^{-}\mathbf{e}_{k}^{+}\right) + \partial_{h}^{-j}g_{i}(mh)f_{k}(mh)\mathbf{e}_{i}^{+}\left(\mathbf{e}_{j}^{-}\mathbf{e}_{k}^{-}\right) \\
\displaystyle \left. + \partial_{h}^{-j}g_{i}(mh)f_{k}(mh)\mathbf{e}_{i}^{-}\left(\mathbf{e}_{j}^{-}\mathbf{e}_{k}^{+}\right) + \partial_{h}^{-j}g_{i}(mh)f_{k}(mh)\mathbf{e}_{i}^{-}\left(\mathbf{e}_{j}^{-}\mathbf{e}_{k}^{-}\right) \right] h^{8},
\end{array}
\end{equation*}
which can be further compressed to the form
\begin{equation*}
\begin{array}{c}
\displaystyle -\sum\limits_{m\in \mathbb{Z}^{8}} \sum\limits_{j=0}^{7}\sum \limits_{i=0}^{7} \sum \limits_{k=0}^{7} \left[\partial_{h}^{+j}g_{i}(mh)f_{k}(mh)\mathbf{e}_{i}\left(\mathbf{e}_{j}^{+}\mathbf{e}_{k}\right) + \partial_{h}^{-j}g_{i}(mh)f_{k}(mh)\mathbf{e}_{i}\left(\mathbf{e}_{j}^{-}\mathbf{e}_{k}\right) \right] h^{8} + 
\\
\displaystyle +2\sum_{m\in \mathbb{Z}^{8}} \sum\limits_{s=1}^{7}  \sum_{i\in I_{s}} \sum_{\stackrel{j\in I_{s}}{j\neq i}}^{7} \sum_{\stackrel{k=1}{k\notin I_{s}}}^{7} \left[\partial_{h}^{+j}g_{i}(mh)f_{k}(mh)\mathbf{e}_{i}\left(\mathbf{e}_{j}^{+}\mathbf{e}_{k}\right) + \partial_{h}^{-j}g_{i}(mh)f_{k}(mh)\mathbf{e}_{i}\left(\mathbf{e}_{j}^{+}\mathbf{e}_{k}\right)\right] h^{8}.
\end{array}
\end{equation*}
Finally, the definition of the Cauchy-Riemann operator $D_{h}^{+-}$ is used:
\begin{equation*}
\begin{array}{c}
\displaystyle -\sum\limits_{m\in \mathbb{Z}^{8}} g(mh) \left[D_{h}^{+-}f(mh)\right] h^{8} + \\
\displaystyle +2\sum_{m\in \mathbb{Z}^{8}} \sum\limits_{s=1}^{7}  \sum_{i\in I_{s}} \sum_{\stackrel{j\in I_{s}}{j\neq i}}^{7} \sum_{\stackrel{k=1}{k\notin I_{s}}}^{7} \left[\partial_{h}^{+j}g_{i}(mh)f_{k}(mh)\mathbf{e}_{i}\left(\mathbf{e}_{j}^{+}\mathbf{e}_{k}\right) + \partial_{h}^{-j}g_{i}(mh)f_{k}(mh)\mathbf{e}_{i}\left(\mathbf{e}_{j}^{+}\mathbf{e}_{k}\right)\right] h^{8}.
\end{array}
\end{equation*}
Thus, the theorem is proved.
\end{proof}\par
Now we will present the discrete Stokes' formulae for the upper half-space $h\mathbb{Z}_{+}^{8}$ and for the lower half-space $h\mathbb{Z}_{-}^{8}$. The proofs of these formulae are generally follow the same strategy as shown above, and the main difference is now that presence of boundary layers, as indicated by Definition~\ref{Definition_geometry}. Thus, we have the following theorem:
\begin{theorem}\label{Theorem:Stokes_half_spaces}
The discrete Stokes' formula for the upper half-space $h\mathbb{Z}_{+}^{8}$ is given by
\begin{equation*}
\begin{array}{c}
\displaystyle \sum_{m\in \mathbb{Z}_{+}^{8}}  \left[ \left( g(mh)D_h^{-+}\right) f(mh) + g(mh) \left( D_h^{+-}f(mh) \right) \right] h^8 =  \\
\displaystyle = 2\sum_{m\in \mathbb{Z}_{+}^{8}} \sum\limits_{s=1}^{7}  \sum_{i\in I_{s}} \sum_{\stackrel{j\in I_{s}}{j\neq i}}^{7} \sum_{\stackrel{k=1}{k\notin I_{s}}}^{7} \left[g_{i}(mh)\mathbf{e}_{i}\left(\partial_{h}^{+j}\mathbf{e}_{j}^{+}f_{k}(mh)\mathbf{e}_{k}\right) \right. \\
\displaystyle \left. +g_{i}(mh)\mathbf{e}_{i}\left(\partial_{h}^{-j}\mathbf{e}_{j}^{-}f_{k}(mh)\mathbf{e}_{k}\right)\right]h^{8} \\
\displaystyle +2\sum\limits_{\underline{m}\in\mathbb{Z}^{7}}\left[\sum_{i=1,6} \sum_{\stackrel{k=1}{k\neq i}}^{6} \left(g_{i}(\underline{m}h,0)\mathbf{e}_{i}\left(\mathbf{e}_{7}^{+}f_{k}(\underline{m}h,h)\mathbf{e}_{k}\right) + g_{i}(\underline{m}h,h)\mathbf{e}_{i}\left(\mathbf{e}_{7}^{-}f_{k}(\underline{m}h,0)\mathbf{e}_{k}\right)\right)\right. \\
\displaystyle \sum_{i=2,5} \sum_{\stackrel{k=1}{k\neq i}}^{6} \left(g_{i}(\underline{m}h,0)\mathbf{e}_{i}\left(\mathbf{e}_{7}^{+}f_{k}(\underline{m}h,h)\mathbf{e}_{k}\right) + g_{i}(\underline{m}h,h)\mathbf{e}_{i}\left(\mathbf{e}_{7}^{-}f_{k}(\underline{m}h,0)\mathbf{e}_{k}\right)\right)\\
\displaystyle \left. \sum_{i=3,4} \sum_{\stackrel{k=1}{k\neq i}}^{6} \left(g_{i}(\underline{m}h,0)\mathbf{e}_{i}\left(\mathbf{e}_{7}^{+}f_{k}(\underline{m}h,h)\mathbf{e}_{k}\right) + g_{i}(\underline{m}h,h)\mathbf{e}_{i}\left(\mathbf{e}_{7}^{-}f_{k}(\underline{m}h,0)\mathbf{e}_{k}\right)\right)\right] h^{7},
\end{array}  
\end{equation*}
and the discrete Stokes' formula for the lower half-space $h\mathbb{Z}_{-}^{8}$ is 
\begin{equation*}
\begin{array}{c}
\displaystyle \sum_{m\in \mathbb{Z}_{-}^{8}}  \left[ \left( g(mh)D_h^{-+}\right) f(mh) + g(mh) \left( D_h^{+-}f(mh) \right) \right] h^8 =  \\
\displaystyle = 2\sum_{m\in \mathbb{Z}_{-}^{8}} \sum\limits_{s=1}^{7}  \sum_{i\in I_{s}} \sum_{\stackrel{j\in I_{s}}{j\neq i}}^{7} \sum_{\stackrel{k=1}{k\notin I_{s}}}^{7} \left[g_{i}(mh)\mathbf{e}_{i}\left(\partial_{h}^{+j}\mathbf{e}_{j}^{+}f_{k}(mh)\mathbf{e}_{k}\right) \right. \\
\displaystyle \left. +g_{i}(mh)\mathbf{e}_{i}\left(\partial_{h}^{-j}\mathbf{e}_{j}^{-}f_{k}(mh)\mathbf{e}_{k}\right)\right]h^{8} \\
\displaystyle +2\sum\limits_{\underline{m}\in\mathbb{Z}^{7}}\left[\sum_{i=1,6} \sum_{\stackrel{k=1}{k\neq i}}^{6} \left(g_{i}(\underline{m}h,-h)\mathbf{e}_{i}\left(\mathbf{e}_{7}^{+}f_{k}(\underline{m}h,0)\mathbf{e}_{k}\right) + g_{i}(\underline{m}h,0)\mathbf{e}_{i}\left(\mathbf{e}_{7}^{-}f_{k}(\underline{m}h,-h)\mathbf{e}_{k}\right)\right)\right. \\
\displaystyle \sum_{i=2,5} \sum_{\stackrel{k=1}{k\neq i}}^{6} \left(g_{i}(\underline{m}h,-h)\mathbf{e}_{i}\left(\mathbf{e}_{7}^{+}f_{k}(\underline{m}h,0)\mathbf{e}_{k}\right) + g_{i}(\underline{m}h,0)\mathbf{e}_{i}\left(\mathbf{e}_{7}^{-}f_{k}(\underline{m}h,-h)\mathbf{e}_{k}\right)\right)\\
\displaystyle \left. \sum_{i=3,4} \sum_{\stackrel{k=1}{k\neq i}}^{6} \left(g_{i}(\underline{m}h,-h)\mathbf{e}_{i}\left(\mathbf{e}_{7}^{+}f_{k}(\underline{m}h,0)\mathbf{e}_{k}\right) + g_{i}(\underline{m}h,0)\mathbf{e}_{i}\left(\mathbf{e}_{7}^{-}f_{k}(\underline{m}h,-h)\mathbf{e}_{k}\right)\right)\right] h^{7},
\end{array}  
\end{equation*}
for all discrete functions $f$ and $g$ such that the series converge. The index sets $I_{s}$, $s=1,\ldots,7$ in these formulae are defined as follows
\begin{equation*}
\begin{array}{cclcclcclccl}
I_{1} & := & \left\{1,2,4\right\}, & I_{2} & := & \left\{1,3,5\right\}, & I_{3} & := & \left\{1,6,7\right\}, & I_{4} & := & \left\{2,3,6\right\}, \\
I_{5} & := & \left\{2,5,7\right\}, & I_{6} & := & \left\{3,4,7\right\}, & I_{7} & := & \left\{4,5,6\right\}.
\end{array}
\end{equation*}
\end{theorem}\par
The discrete Stokes' formulae presented in Theorem~\ref{Theorem:Stokes_half_spaces} naturally allows us to introduce the discrete octonionic Borel-Pompeiu formulae for the lower and upper half-spaces. This is done by substituting the function $g$ in these formulae via the shifted discrete fundamental solution $E_{h}^{-+}(\cdot - mh)$ with $m\in\mathbb{Z}_{\pm}^{8}$ and then taking into account Definition~\ref{Definition:Discrete_fundamental_solution}. After that, by considering discrete octonionic monogenic functions, discrete octonionic Cauchy formulae can be straightforwardly defined. Because the machinery is analogous, we will only present here the discrete octonionic Borel-Pompeiu and Cauchy formulae for the upper half-space:
\begin{theorem}
Let $E_{h}^{-+}$ be the discrete fundamental solution of the discrete Cauchy-Riemann operator $D_{h}^{-+}$. Then the discrete octonionic Borel-Pompeiu formula for the upper half-lattice $h\mathbb{Z}_{+}^{8}$ is given by
\begin{equation*}
\begin{array}{c}
\displaystyle \sum_{n\in \mathbb{Z}_{+}^{8}}  E_{h}^{-+}(nh-mh) \left[ D_h^{+-}f(nh) \right] h^8 \\
\displaystyle -2\sum_{m\in \mathbb{Z}_{+}^{8}} \sum\limits_{s=1}^{7}  \sum_{i\in I_{s}} \sum_{\stackrel{j\in I_{s}}{j\neq i}}^{7} \sum_{\stackrel{k=1}{k\notin I_{s}}}^{7} \left[E_{h,i}^{-+}(nh-mh)(mh)\mathbf{e}_{i}\left(\partial_{h}^{+j}\mathbf{e}_{j}^{+}f_{k}(nh)\mathbf{e}_{k}\right) \right. \\
\displaystyle \left. +E_{h,i}^{-+}(nh-mh)\mathbf{e}_{i}\left(\partial_{h}^{-j}\mathbf{e}_{j}^{-}f_{k}(nh)\mathbf{e}_{k}\right)\right]h^{8} \\
\displaystyle -2\sum\limits_{\underline{m}\in\mathbb{Z}^{7}}\left[\sum_{i=1,6} \sum_{\stackrel{k=1}{k\neq i}}^{6} \left(E_{h,i}^{-+}(\underline{n}h-\underline{m}h,0)\mathbf{e}_{i}\left(\mathbf{e}_{7}^{+}f_{k}(\underline{n}h,h)\mathbf{e}_{k}\right) \right. \right. \\
\displaystyle \left. + E_{h,i}^{-+}(\underline{n}h-\underline{m}h,h)\mathbf{e}_{i}\left(\mathbf{e}_{7}^{-}f_{k}(\underline{n}h,0)\mathbf{e}_{k}\right)\right) \\
\displaystyle \sum_{i=2,5} \sum_{\stackrel{k=1}{k\neq i}}^{6} \left(E_{h,i}^{-+}(\underline{n}h-\underline{m}h,0)\mathbf{e}_{i}\left(\mathbf{e}_{7}^{+}f_{k}(\underline{n}h,h)\mathbf{e}_{k}\right) \right. \\
\displaystyle \left. + E_{h,i}^{-+}(\underline{n}h-\underline{m}h,h)\mathbf{e}_{i}\left(\mathbf{e}_{7}^{-}f_{k}(\underline{n}h,0)\mathbf{e}_{k}\right)\right)\\
\displaystyle \sum_{i=3,4} \sum_{\stackrel{k=1}{k\neq i}}^{6} \left(E_{h,i}^{-+}(\underline{n}h-\underline{m}h,0)\mathbf{e}_{i}\left(\mathbf{e}_{7}^{+}f_{k}(\underline{m}h,h)\mathbf{e}_{k}\right) \right. \\
\displaystyle \left. \left. + E_{h,i}^{-+}(\underline{n}h-\underline{m}h,h)\mathbf{e}_{i}\left(\mathbf{e}_{7}^{-}f_{k}(\underline{m}h,0)\mathbf{e}_{k}\right)\right)\right] h^{7} = \left\{\begin{array}{cl}
0, & m\notin\mathbb{Z}_{+}^{8}, \\
-f(mh), & m\in\mathbb{Z}_{+}^{8},
\end{array} \right.
\end{array}
\end{equation*}
for any discrete function $f$ such that the series converge, and index sets $I_{s}$, $s=1,\ldots,7$ are given by
\begin{equation*}
\begin{array}{cclcclcclccl}
I_{1} & := & \left\{1,2,4\right\}, & I_{2} & := & \left\{1,3,5\right\}, & I_{3} & := & \left\{1,6,7\right\}, & I_{4} & := & \left\{2,3,6\right\}, \\
I_{5} & := & \left\{2,5,7\right\}, & I_{6} & := & \left\{3,4,7\right\}, & I_{7} & := & \left\{4,5,6\right\}.
\end{array}
\end{equation*} 
In the case when $f$ is a discrete octonionic left monogenic function with respect to the operator $D_{h}^{+-}$, we obtain the discrete octonionic Cauchy formula for the upper half-lattice $h\mathbb{Z}_{+}^{8}$ in the form 
\begin{equation*}
\begin{array}{c}
\displaystyle 2\sum_{m\in \mathbb{Z}_{+}^{8}} \sum\limits_{s=1}^{7}  \sum_{i\in I_{s}} \sum_{\stackrel{j\in I_{s}}{j\neq i}}^{7} \sum_{\stackrel{k=1}{k\notin I_{s}}}^{7} \left[E_{h,i}^{-+}(nh-mh)(mh)\mathbf{e}_{i}\left(\partial_{h}^{+j}\mathbf{e}_{j}^{+}f_{k}(nh)\mathbf{e}_{k}\right) \right. \\
\displaystyle \left. +E_{h,i}^{-+}(nh-mh)\mathbf{e}_{i}\left(\partial_{h}^{-j}\mathbf{e}_{j}^{-}f_{k}(nh)\mathbf{e}_{k}\right)\right]h^{8} \\
\displaystyle 2\sum\limits_{\underline{m}\in\mathbb{Z}^{7}}\left[\sum_{i=1,6} \sum_{\stackrel{k=1}{k\neq i}}^{6} \left(E_{h,i}^{-+}(\underline{n}h-\underline{m}h,0)\mathbf{e}_{i}\left(\mathbf{e}_{7}^{+}f_{k}(\underline{n}h,h)\mathbf{e}_{k}\right) \right. \right. \\
\displaystyle \left. + E_{h,i}^{-+}(\underline{n}h-\underline{m}h,h)\mathbf{e}_{i}\left(\mathbf{e}_{7}^{-}f_{k}(\underline{n}h,0)\mathbf{e}_{k}\right)\right) \\
\displaystyle \sum_{i=2,5} \sum_{\stackrel{k=1}{k\neq i}}^{6} \left(E_{h,i}^{-+}(\underline{n}h-\underline{m}h,0)\mathbf{e}_{i}\left(\mathbf{e}_{7}^{+}f_{k}(\underline{n}h,h)\mathbf{e}_{k}\right) \right. \\
\displaystyle \left. + E_{h,i}^{-+}(\underline{n}h-\underline{m}h,h)\mathbf{e}_{i}\left(\mathbf{e}_{7}^{-}f_{k}(\underline{n}h,0)\mathbf{e}_{k}\right)\right)\\
\displaystyle \sum_{i=3,4} \sum_{\stackrel{k=1}{k\neq i}}^{6} \left(E_{h,i}^{-+}(\underline{n}h-\underline{m}h,0)\mathbf{e}_{i}\left(\mathbf{e}_{7}^{+}f_{k}(\underline{m}h,h)\mathbf{e}_{k}\right) \right. \\
\displaystyle \left. \left. + E_{h,i}^{-+}(\underline{n}h-\underline{m}h,h)\mathbf{e}_{i}\left(\mathbf{e}_{7}^{-}f_{k}(\underline{m}h,0)\mathbf{e}_{k}\right)\right)\right] h^{7} = \left\{\begin{array}{cl}
0, & m\notin\mathbb{Z}_{+}^{8}, \\
f(mh), & m\in\mathbb{Z}_{+}^{8}.
\end{array} \right.
\end{array}
\end{equation*}
\end{theorem}\par
The discrete octonionic Cauchy formula naturally provides a definition of a Cauchy transform, which we will not write explicitly here, but refer to \cite{CKKS} for the case of discrete Clifford analysis on how to proceed. After that discrete Hardy spaces for half-spaces can be introduced. We will discuss the construction of discrete Hardy spaces for bounded domains in Section~\ref{Section:Hardy_bounded}, where we will also briefly comment the case of half-spaces.\par
It is also important to underline that while considering the discrete formulae for half-spaces, a three-layer structure of the discrete boundary appears naturally in the constructions, as introduced in Definition~\ref{Definition_geometry}. Additionally, the associator term appears also in all constructions, as in the classical continuous case.\par

\section{Discrete octonionic Stokes' formula for bounded domains}\label{Section:Stokes}

In this section, we will present the discrete octonionic  Stokes' formula for bounded domains. As it has been mentioned previously, there are generally two approaches to construct the discrete Stokes' formula, namely either by the help of explicit constructions, or more abstractly by using characteristic functions. Since our aim here is first to study the influence of the non-associativity of the octonionic multiplication, we will first construct the discrete Stokes' formula by the constructive approach for a cuboid in $\mathbb{R}^{8}$, see also \cite{Cerejeiras_2}. After that, we will extend this result to arbitrary bounded domains in $\mathbb{R}^{8}$ by the help of characteristic functions.\par
Let now $K_{h}$ be a uniformly meshed cuboid in $\mathbb{R}^{8}$. Moreover, we assume the following index convention for the interior points of $K_{h}$:
\begin{equation*}
K_{h} := \left\{(m_{0},m_{1},\ldots,m_{7}) \mid  m_{i}\in[1,N_{i}-1], m_{i}\in\mathbb{N}, i=0,1,\ldots,7  \right\}.
\end{equation*}
Hence, the points with indices $m_{i}=0$ and $m_{i}=N_{i}$, $i=0,1,\ldots,7$ represent boundary layer $\gamma_{h}^{*}$. Moreover, for convenience reason, we sub-divide the boundary layer $\gamma_{h}^{*}$ into parts being parallel to the coordinate axes, denoted as follows:
\begin{equation*}
\gamma_{h,2i}^{*}, \quad \gamma_{h,2i+1}^{*}, \quad i=0,1,\ldots,7.
\end{equation*}\par
Now, we present the following theorem:
\begin{theorem}\label{Theorem:Stokes_cuboid}
The discrete Stokes' formula for the cuboid $K_{h}$ with a uniform lattice is given by
\begin{equation}
\label{DiscreteStokesFormula_cuboid}
\begin{array}{c}
\displaystyle \sum\limits_{m\in K_{h}}  \left[ \left( g(mh)D_h^{-+}\right) f(mh) + g(mh) \left( D_h^{+-}f(mh) \right)  \right] h^{8} = \\
\displaystyle = -\sum\limits_{mh\in\gamma_{h,2j}^{*}}\sum\limits_{j=0}^{7} \left[g(mh)\mathbf{e}_{j}^{+}f(mh+\mathbf{e}_{j}h) + g(mh+\mathbf{e}_{j}h)\mathbf{e}_{j}^{-}f(mh) \right]h^{7} \\
\displaystyle + \sum\limits_{mh\in\gamma_{h,2j+1}^{*}}\sum\limits_{j=0}^{7} \left[g(mh-\mathbf{e}_{j}h)\mathbf{e}_{j}^{+}f(mh) + g(mh)\mathbf{e}_{j}^{-}f(mh-\mathbf{e}_{j}h) \right]h^{7} \\
\displaystyle + 2\sum_{m\in K_{h}} \sum\limits_{s=1}^{7}  \sum_{i\in I_{s}} \sum_{\stackrel{j\in I_{s}}{j\neq i}}^{7} \sum_{\stackrel{k=1}{k\notin I_{s}}}^{7} \left[g_{i}(mh)\mathbf{e}_{i}\left(\partial_{h}^{+j}\mathbf{e}_{j}^{+}f_{k}(mh)\mathbf{e}_{k}\right) \right. \\
\displaystyle \left. + g_{i}(mh)\mathbf{e}_{i}\left(\partial_{h}^{-j}\mathbf{e}_{j}^{-}f_{k}(mh)\mathbf{e}_{k}\right)\right]h^{8} \\
\displaystyle + 2\sum\limits_{mh\in\gamma_{h,2j}^{*}} \sum\limits_{s=1}^{7}  \sum_{i\in I_{s}} \sum_{\stackrel{j\in I_{s}}{j\neq i}}^{7} \sum_{\stackrel{k=1}{k\notin I_{s}}}^{7} \left[g_{i}(mh)\mathbf{e}_{i}\left(\mathbf{e}_{j}^{+}f_{k}(mh+\mathbf{e}_{j}h)\mathbf{e}_{k}\right) \right. \\
\displaystyle \left. + g_{i}(mh+\mathbf{e}_{j}h)\mathbf{e}_{i}\left(\mathbf{e}_{j}^{-}f_{k}(mh)\mathbf{e}_{k}\right)\right]h^{7} \\
\displaystyle - 2\sum\limits_{mh\in\gamma_{h,2j+1}^{*}} \sum\limits_{s=1}^{7}  \sum_{i\in I_{s}} \sum_{\stackrel{j\in I_{s}}{j\neq i}}^{7} \sum_{\stackrel{k=1}{k\notin I_{s}}}^{7} \left[g_{i}(mh-\mathbf{e}_{j}h)\mathbf{e}_{i}\left(\mathbf{e}_{j}^{+}f_{k}(mh)\mathbf{e}_{k}\right) \right. \\
\displaystyle \left. + g_{i}(mh)\mathbf{e}_{i}\left(\mathbf{e}_{j}^{-}f_{k}(mh-\mathbf{e}_{j}h)\mathbf{e}_{k}\right)\right]h^{7},
\end{array}  
\end{equation}
for all discrete functions $f$ and $g$ such that the series converge, where the last three sums are related to the associator, and with $\gamma_{h,2j}^{*}$ and $\gamma_{h,2j+1}^{*}$ being \textquotedblleft left\textquotedblright\, and \textquotedblleft right\textquotedblright\, boundary parts of $K_{h}$, respectively. The index sets $I_{s}$, $s=1,\ldots,7$ in this formula are defined as follows
\begin{equation*}
\begin{array}{cclcclcclccl}
I_{1} & := & \left\{1,2,4\right\}, & I_{2} & := & \left\{1,3,5\right\}, & I_{3} & := & \left\{1,6,7\right\}, & I_{4} & := & \left\{2,3,6\right\}, \\
I_{5} & := & \left\{2,5,7\right\}, & I_{6} & := & \left\{3,4,7\right\}, & I_{7} & := & \left\{4,5,6\right\}.
\end{array}
\end{equation*}
\end{theorem}
\begin{proof}
The proof of the discrete Stokes' formula for a bounded cuboid $K_{h}$ follows the steps of the proof presented for Theorem~\ref{Discrete_Stokes_space}. The main difficulty is to address the different boundary parts. Hence, we will present only essential steps here. The proof starts again by working with the first term on the left-hand side of~(\ref{DiscreteStokesFormula_cuboid}). After some simplifications and using the same trick as in the proof of Theorem~\ref{Discrete_Stokes_space} by adding and subtracting the non-associative terms, we obtain the following expression
\begin{equation*}
\begin{array}{c}
\displaystyle \sum\limits_{m\in K_{h}} \sum\limits_{j=0}^{7}\sum \limits_{i=0}^{7} \sum \limits_{k=0}^{7} \left[\partial_{h}^{-j}g_{i}(mh)f_{k}(mh)\mathbf{e}_{i}^{+}\left(\mathbf{e}_{j}^{+}\mathbf{e}_{k}^{+}\right) + \partial_{h}^{-j}g_{i}(mh)f_{k}(mh)\mathbf{e}_{i}^{+}\left(\mathbf{e}_{j}^{+}\mathbf{e}_{k}^{-}\right) \right. \\
\displaystyle + \partial_{h}^{-j}g_{i}(mh)f_{k}(mh)\mathbf{e}_{i}^{-}\left(\mathbf{e}_{j}^{+}\mathbf{e}_{k}^{+}\right) + \partial_{h}^{-j}g_{i}(mh)f_{k}(mh)\mathbf{e}_{i}^{-}\left(\mathbf{e}_{j}^{+}\mathbf{e}_{k}^{-}\right) \\
\displaystyle + \partial_{h}^{+j}g_{i}(mh)f_{k}(mh)\mathbf{e}_{i}^{+}\left(\mathbf{e}_{j}^{-}\mathbf{e}_{k}^{+}\right) + \partial_{h}^{+j}g_{i}(mh)f_{k}(mh)\mathbf{e}_{i}^{+}\left(\mathbf{e}_{j}^{-}\mathbf{e}_{k}^{-}\right) \\
\displaystyle \left. + \partial_{h}^{+j}g_{i}(mh)f_{k}(mh)\mathbf{e}_{i}^{-}\left(\mathbf{e}_{j}^{-}\mathbf{e}_{k}^{+}\right) + \partial_{h}^{+j}g_{i}(mh)f_{k}(mh)\mathbf{e}_{i}^{-}\left(\mathbf{e}_{j}^{-}\mathbf{e}_{k}^{-}\right) \right] \\
\displaystyle - 2\sum_{m\in K_{h}} \sum\limits_{s=1}^{7}  \sum_{i\in I_{s}} \sum_{\stackrel{j\in I_{s}}{j\neq i}}^{7} \sum_{\stackrel{k=1}{k\notin I_{s}}}^{7} \left[\partial_{h}^{-j}g_{i}(mh)f_{k}(mh)\mathbf{e}_{i}^{+}\left(\mathbf{e}_{j}^{+}\mathbf{e}_{k}^{+}\right) + \partial_{h}^{-j}g_{i}(mh)f_{k}(mh)\mathbf{e}_{i}^{+}\left(\mathbf{e}_{j}^{+}\mathbf{e}_{k}^{-}\right) \right. \\
\displaystyle + \partial_{h}^{-j}g_{i}(mh)f_{k}(mh)\mathbf{e}_{i}^{-}\left(\mathbf{e}_{j}^{+}\mathbf{e}_{k}^{+}\right) + \partial_{h}^{-j}g_{i}(mh)f_{k}(mh)\mathbf{e}_{i}^{-}\left(\mathbf{e}_{j}^{+}\mathbf{e}_{k}^{-}\right) \\
\displaystyle + \partial_{h}^{+j}g_{i}(mh)f_{k}(mh)\mathbf{e}_{i}^{+}\left(\mathbf{e}_{j}^{-}\mathbf{e}_{k}^{+}\right) + \partial_{h}^{+j}g_{i}(mh)f_{k}(mh)\mathbf{e}_{i}^{+}\left(\mathbf{e}_{j}^{-}\mathbf{e}_{k}^{-}\right) \\
\displaystyle \left. + \partial_{h}^{+j}g_{i}(mh)f_{k}(mh)\mathbf{e}_{i}^{-}\left(\mathbf{e}_{j}^{-}\mathbf{e}_{k}^{+}\right) + \partial_{h}^{+j}g_{i}(mh)f_{k}(mh)\mathbf{e}_{i}^{-}\left(\mathbf{e}_{j}^{-}\mathbf{e}_{k}^{-}\right)\right]h^{8},
\end{array}
\end{equation*}
where the index sets $I_{s},s=1,\ldots,7$ have been already introduced before.\par
The next step is to use the definition of the finite difference operators and perform the change of variables. The critical point here is keep in mind, that $K_{h}$, according to our convention, contains only points with indices $m_{i}\in[1,N_{i}-1]$, $i=0,1,\ldots,7$. This fact implies that boundary points will appear after using the finite difference operators and performing change of variables. By using the same ideas as in the proof of Theorem~\ref{Discrete_Stokes_space}, the summation over $K_{h}$ can be written in terms of the operator $D_{h}^{+-}$, which is not possible to do for the summation over boundary points and for the associator terms. Thus, we obtain the following expression:
\begin{equation*}
\begin{array}{c}
\displaystyle -\sum\limits_{m\in K_{h}} g(mh)\left[D_{h}^{+-}f(mh)\right]h^{8} \\
-\displaystyle \sum\limits_{m\in K_{h}} \sum\limits_{j=0}^{7}\sum \limits_{i=0}^{7} \sum \limits_{k=0}^{7} \left[ -g_{i}(0_{j})\mathbf{e}_{i}f_{k}(1_{j})\left(\mathbf{e}_{j}^{+}\mathbf{e}_{k}\right) - g_{i}(N_{j}-1)\mathbf{e}_{i}f_{k}(N_{j})\left(\mathbf{e}_{j}^{+}\mathbf{e}_{k}\right)\right. \\
\displaystyle \left. - g_{i}(N_{j})\mathbf{e}_{i}f_{k}(N_{j}-1)\left(\mathbf{e}_{j}^{-}\mathbf{e}_{k}\right) + g_{i}(1_{j})\mathbf{e}_{i}f_{k}(0_{j})\left(\mathbf{e}_{j}^{-}\mathbf{e}_{k}\right)\right]h^{7} \\
\displaystyle +2\sum_{m\in K_{h}} \sum\limits_{s=1}^{7}  \sum_{i\in I_{s}} \sum_{\stackrel{j\in I_{s}}{j\neq i}}^{7} \sum_{\stackrel{k=1}{k\notin I_{s}}}^{7} \left[g_{i}(mh)\mathbf{e}_{i}\partial_{h}^{+j}f_{k}(mh)\left(\mathbf{e}_{j}^{+}\mathbf{e}_{k}\right) + g_{i}(mh)\mathbf{e}_{i}\partial_{h}^{-j}f_{k}(mh)\left(\mathbf{e}_{j}^{-}\mathbf{e}_{k}\right) \right]h^{8} \\
\displaystyle \displaystyle +2\sum_{m\in K_{h}} \sum\limits_{s=1}^{7}  \sum_{i\in I_{s}} \sum_{\stackrel{j\in I_{s}}{j\neq i}}^{7} \sum_{\stackrel{k=1}{k\notin I_{s}}}^{7} \left[g_{i}(0_{j})\mathbf{e}_{i}f_{k}(1_{j})\left(\mathbf{e}_{j}^{+}\mathbf{e}_{k}\right) - g_{i}(N_{j}-1)\mathbf{e}_{i}f_{k}(N_{j})\left(\mathbf{e}_{j}^{+}\mathbf{e}_{k}\right) \right. \\
\displaystyle \left. g_{i}(N_{j})\mathbf{e}_{i}f_{k}(N_{j}-1)\left(\mathbf{e}_{j}^{-}\mathbf{e}_{k}\right) + g_{i}(1_{j})\mathbf{e}_{i}f_{k}(0_{j})\left(\mathbf{e}_{j}^{-}\mathbf{e}_{k}\right) \right]h^{7},
\end{array}
\end{equation*}
where the sub-index $j$ by $0,1,N_{j}$ indicates which coordinate $m_{j}$ must be set to $0,1$ or $N$, respectively. For example, $g(0_{3})$ implies the following full notation
\begin{equation*}
g(m_{0}h,m_{1}h,m_{2}h,0,m_{4}h,m_{5}h,m_{6}h,m_{7}h).
\end{equation*}
Further, it is important to remark, that the summations which include $g_{i}(0)$, $g_{i}(1)$, $g_{i}(N_{j}-1)$, $g_{i}(N_{j})$, $f_{k}(0)$, $f_{k}(1)$, $f_{k}(N_{j}-1)$, $f_{k}(N_{j})$ represent summations over different boundary parts. These summations can be written by the help of boundary layers $\gamma_{h,2i}^{*}$ and $\gamma_{h,2i+1}^{+}$. By writing these sums, we obtain the discrete octonionic Stokes' formula for a bounded cuboid $K_{h}$. Thus, the theorem is proved.
\end{proof}\par
As the next step, we present an extension of Theorem~\ref{Theorem:Stokes_cuboid} to arbitrary bounded domains in $\mathbb{R}^{8}$ by the help of characteristic functions and ideas presented in \cite{Cerejeiras}. This extension is straightforwardly possible, because of working with boundary parts $\gamma_{h,2j}$ and $\gamma_{h,2j+1}$, which do not impose any particular indexing strategy in comparison to the construction of $K_{h}$. Hence, we have the following theorem: 
\begin{theorem}\label{Theorem:Stokes_interior_arbitrary}
Let $\Omega_{h} \subset h\mathbb{Z}^{8}$ be a discrete simply connected bounded domain, then the following formula holds
\begin{equation}
\label{DiscreteStokesFormula_interior_arbitrary}
\begin{array}{c}
\displaystyle \sum\limits_{m\in \mathbb{Z}^{8}}  \left[ \left( g(mh)D_h^{-+}\right) f(mh) + g(mh) \left( D_h^{+-}f(mh) \right)  \right]\chi_{\Omega_{h}} h^{8} = \\
\displaystyle = -\sum\limits_{mh\in\gamma_{h,2j}^{*}}\sum\limits_{j=0}^{7} \left[g(mh)\mathbf{e}_{j}^{+}f(mh+\mathbf{e}_{j}h) + g(mh+\mathbf{e}_{j}h)\mathbf{e}_{j}^{-}f(mh) \right]h^{7} \\
\displaystyle + \sum\limits_{mh\in\gamma_{h,2j+1}^{*}}\sum\limits_{j=0}^{7} \left[g(mh-\mathbf{e}_{j}h)\mathbf{e}_{j}^{+}f(mh) + g(mh)\mathbf{e}_{j}^{-}f(mh-\mathbf{e}_{j}h) \right]h^{7} \\
\displaystyle + 2\sum_{m\in \mathbb{Z}^{8}} \sum\limits_{s=1}^{7}  \sum_{i\in I_{s}} \sum_{\stackrel{j\in I_{s}}{j\neq i}}^{7} \sum_{\stackrel{k=1}{k\notin I_{s}}}^{7} \left[g_{i}(mh)\mathbf{e}_{i}\left(\partial_{h}^{+j}\mathbf{e}_{j}^{+}f_{k}(mh)\mathbf{e}_{k}\right) \right. \\
\displaystyle \left. + g_{i}(mh)\mathbf{e}_{i}\left(\partial_{h}^{-j}\mathbf{e}_{j}^{-}f_{k}(mh)\mathbf{e}_{k}\right)\right]\chi_{\Omega_{h}}h^{8} \\
\displaystyle + 2\sum\limits_{mh\in\gamma_{h,2j}^{*}} \sum\limits_{s=1}^{7}  \sum_{i\in I_{s}} \sum_{\stackrel{j\in I_{s}}{j\neq i}}^{7} \sum_{\stackrel{k=1}{k\notin I_{s}}}^{7} \left[g_{i}(mh)\mathbf{e}_{i}\left(\mathbf{e}_{j}^{+}f_{k}(mh+\mathbf{e}_{j}h)\mathbf{e}_{k}\right) \right. \\
\displaystyle \left. + g_{i}(mh+\mathbf{e}_{j}h)\mathbf{e}_{i}\left(\mathbf{e}_{j}^{-}f_{k}(mh)\mathbf{e}_{k}\right)\right]h^{7} \\
\displaystyle - 2\sum\limits_{mh\in\gamma_{h,2j+1}^{*}} \sum\limits_{s=1}^{7}  \sum_{i\in I_{s}} \sum_{\stackrel{j\in I_{s}}{j\neq i}}^{7} \sum_{\stackrel{k=1}{k\notin I_{s}}}^{7} \left[g_{i}(mh-\mathbf{e}_{j}h)\mathbf{e}_{i}\left(\mathbf{e}_{j}^{+}f_{k}(mh)\mathbf{e}_{k}\right) \right. \\
\displaystyle \left. + g_{i}(mh)\mathbf{e}_{i}\left(\mathbf{e}_{j}^{-}f_{k}(mh-\mathbf{e}_{j}h)\mathbf{e}_{k}\right)\right]h^{7},
\end{array}
\end{equation}
for all discrete functions $f$ and $g$ such that the series converge, and where $\chi_{\Omega_{h}}$ is the characteristic function of the discrete domain, and the index sets are given by
\begin{equation*}
\begin{array}{cclcclcclccl}
I_{1} & := & \left\{1,2,4\right\}, & I_{2} & := & \left\{1,3,5\right\}, & I_{3} & := & \left\{1,6,7\right\}, & I_{4} & := & \left\{2,3,6\right\}, \\
I_{5} & := & \left\{2,5,7\right\}, & I_{6} & := & \left\{3,4,7\right\}, & I_{7} & := & \left\{4,5,6\right\}.
\end{array}
\end{equation*}
\end{theorem}\par
Finally, a similar formula can be obtained for a discrete exterior domain $\Omega_{h}^{ext}$, where we need to take into account the reversion of directions:
\begin{theorem}
Let $\Omega_{h}^{ext}$ be the discrete exterior domain associated with $\Omega_{h}$ (see Theorem~\ref{Theorem:Stokes_interior_arbitrary}), then the following formula holds
\begin{equation}
\label{DiscreteStokesFormula_exterior_arbitrary}
\begin{array}{c}
\displaystyle \sum\limits_{m\in \mathbb{Z}^{8}}  \left[ \left( g(mh)D_h^{-+}\right) f(mh) + g(mh) \left( D_h^{+-}f(mh) \right)  \right]\chi_{\Omega_{h}^{ext}} h^{8} = \\
\displaystyle = -\sum\limits_{mh\in\gamma_{h,2j+1}^{*}}\sum\limits_{j=0}^{7} \left[g(mh)\mathbf{e}_{j}^{+}f(mh+\mathbf{e}_{j}h) + g(mh+\mathbf{e}_{j}h)\mathbf{e}_{j}^{-}f(mh) \right]h^{7} \\
\displaystyle + \sum\limits_{mh\in\gamma_{h,2j}^{*}}\sum\limits_{j=0}^{7} \left[g(mh-\mathbf{e}_{j}h)\mathbf{e}_{j}^{+}f(mh) + g(mh)\mathbf{e}_{j}^{-}f(mh-\mathbf{e}_{j}h) \right]h^{7} \\
\displaystyle + 2\sum_{m\in \mathbb{Z}^{8}} \sum\limits_{s=1}^{7}  \sum_{i\in I_{s}} \sum_{\stackrel{j\in I_{s}}{j\neq i}}^{7} \sum_{\stackrel{k=1}{k\notin I_{s}}}^{7} \left[g_{i}(mh)\mathbf{e}_{i}\left(\partial_{h}^{+j}\mathbf{e}_{j}^{+}f_{k}(mh)\mathbf{e}_{k}\right) \right. \\
\displaystyle \left. + g_{i}(mh)\mathbf{e}_{i}\left(\partial_{h}^{-j}\mathbf{e}_{j}^{-}f_{k}(mh)\mathbf{e}_{k}\right)\right]\chi_{\Omega_{h}^{ext}}h^{8} \\
\displaystyle + 2\sum\limits_{mh\in\gamma_{h,2j+1}^{*}} \sum\limits_{s=1}^{7}  \sum_{i\in I_{s}} \sum_{\stackrel{j\in I_{s}}{j\neq i}}^{7} \sum_{\stackrel{k=1}{k\notin I_{s}}}^{7} \left[g_{i}(mh)\mathbf{e}_{i}\left(\mathbf{e}_{j}^{+}f_{k}(mh+\mathbf{e}_{j}h)\mathbf{e}_{k}\right) \right. \\
\displaystyle \left. + g_{i}(mh+\mathbf{e}_{j}h)\mathbf{e}_{i}\left(\mathbf{e}_{j}^{-}f_{k}(mh)\mathbf{e}_{k}\right)\right]h^{7} \\
\displaystyle - 2\sum\limits_{mh\in\gamma_{h,2j}^{*}} \sum\limits_{s=1}^{7}  \sum_{i\in I_{s}} \sum_{\stackrel{j\in I_{s}}{j\neq i}}^{7} \sum_{\stackrel{k=1}{k\notin I_{s}}}^{7} \left[g_{i}(mh-\mathbf{e}_{j}h)\mathbf{e}_{i}\left(\mathbf{e}_{j}^{+}f_{k}(mh)\mathbf{e}_{k}\right) \right. \\
\displaystyle \left. + g_{i}(mh)\mathbf{e}_{i}\left(\mathbf{e}_{j}^{-}f_{k}(mh-\mathbf{e}_{j}h)\mathbf{e}_{k}\right)\right]h^{7},
\end{array}
\end{equation}
for all discrete functions $f$ and $g$ such that the series converge, and where $\chi_{\Omega_{h}^{ext}}$ is the characteristic function of the discrete domain, and the index sets are given by
\begin{equation*}
\begin{array}{cclcclcclccl}
I_{1} & := & \left\{1,2,4\right\}, & I_{2} & := & \left\{1,3,5\right\}, & I_{3} & := & \left\{1,6,7\right\}, & I_{4} & := & \left\{2,3,6\right\}, \\
I_{5} & := & \left\{2,5,7\right\}, & I_{6} & := & \left\{3,4,7\right\}, & I_{7} & := & \left\{4,5,6\right\}.
\end{array}
\end{equation*}
\end{theorem}\par
\begin{remark}
It is important to remark, that the use of characteristic functions is nonetheless strongly based on the constructive approach discussed at first. In fact, characteristic functions provide a convenient way to extend results from simple geometries (cuboid), which are obtained constructively, to arbitrary domains. Further, evidently the non-associativity of octonionic multiplication gives rise to the associator terms, which appear for all interior points, as well as for boundary points. Moreover, as it can be clearly seen from the results presented above, the expressions for the associator terms have much more complicated structures, because several indices are always missing from the summations.
\end{remark}\par

\section{Octonionic Borel-Pompeiu and Cauchy formulae for bounded domains}\label{Section:Borel_Cauchy}

As well-known, a key ingredient of a function theory is the Borel-Pompeiu formula and the Cauchy formula representing its specific case, when a function $f$ is holomorphic (monogenic). Therefore, we introduce these formulae for discrete interior and exterior cases within the setting of bounded domains in $\mathbb{R}^{8}$. The proof of the discrete octonionic Borel-Pompeiu formulae is based on the discrete Stokes' formulae introduced in Section~\ref{Section:Stokes}, where the function $g$ is replaced by the shifted discrete fundamental solution $E_{h}^{-+}(\cdot - mh)$, see also \cite{CKKS,Cerejeiras}. The following theorem present the discrete octonionic Borel-Pompeiu formulae for bounded domain in the exterior and interior settings:
\begin{theorem}
Let $\Omega_{h}\subset h\mathbb{Z}^{8}$ be a discrete simply connected bounded domain and let $\Omega_{h}^{ext}$ be its associated exterior domain. Then the discrete octonionic Borel-Pompeiu formula for $\Omega_{h}$ is given by
\begin{equation*}
\begin{array}{c}
\displaystyle \sum\limits_{r\in \mathbb{Z}^{8}}  E_{h}^{-+}(rh-mh)\left( D_h^{+-}f(rh) \right) \chi_{\Omega_{h}}(rh) h^{8} \\
\displaystyle +\sum\limits_{rh\in\gamma_{h,2j}^{*}}\sum\limits_{j=0}^{7} \left[E_{h}^{-+}(rh-mh)\mathbf{e}_{j}^{+}f(rh+\mathbf{e}_{j}h) + E_{h}^{-+}(rh-mh-\mathbf{e}_{j}h)\mathbf{e}_{j}^{-}f(rh) \right]h^{7} \\
\displaystyle - \sum\limits_{rh\in\gamma_{h,2j+1}^{*}}\sum\limits_{j=0}^{7} \left[E_{h}^{-+}(rh-mh+\mathbf{e}_{j}h)\mathbf{e}_{j}^{+}f(rh) + E_{h}^{-+}(rh-mh)\mathbf{e}_{j}^{-}f(rh-\mathbf{e}_{j}h) \right]h^{7} \\
\displaystyle - 2\sum_{r\in \mathbb{Z}^{8}} \sum\limits_{s=1}^{7}  \sum_{i\in I_{s}} \sum_{\stackrel{j\in I_{s}}{j\neq i}}^{7} \sum_{\stackrel{k=1}{k\notin I_{s}}}^{7} \left[E_{h,i}^{-+}(rh-mh)\mathbf{e}_{i}\left(\partial_{h}^{+j}\mathbf{e}_{j}^{+}f_{k}(rh)\mathbf{e}_{k}\right) \right. \\
\displaystyle \left. + E_{h,i}^{-+}(rh-mh)\mathbf{e}_{i}\left(\partial_{h}^{-j}\mathbf{e}_{j}^{-}f_{k}(rh)\mathbf{e}_{k}\right)\right]\chi_{\Omega_{h}}(rh) h^{8} \\
\displaystyle - 2\sum\limits_{rh\in\gamma_{h,2j}^{*}} \sum\limits_{s=1}^{7}  \sum_{i\in I_{s}} \sum_{\stackrel{j\in I_{s}}{j\neq i}}^{7} \sum_{\stackrel{k=1}{k\notin I_{s}}}^{7} \left[E_{h,i}^{-+}(rh-mh)\mathbf{e}_{i}\left(\mathbf{e}_{j}^{+}f_{k}(rh+\mathbf{e}_{j}h)\mathbf{e}_{k}\right) \right. \\
\displaystyle \left. + E_{h,i}^{-+}(rh-mh-\mathbf{e}_{j}h)\mathbf{e}_{i}\left(\mathbf{e}_{j}^{-}f_{k}(rh)\mathbf{e}_{k}\right)\right]h^{7} \\
\displaystyle + 2\sum\limits_{rh\in\gamma_{h,2j+1}^{*}} \sum\limits_{s=1}^{7}  \sum_{i\in I_{s}} \sum_{\stackrel{j\in I_{s}}{j\neq i}}^{7} \sum_{\stackrel{k=1}{k\notin I_{s}}}^{7} \left[E_{h,i}^{-+}(rh-mh+\mathbf{e}_{j}h)\mathbf{e}_{i}\left(\mathbf{e}_{j}^{+}f_{k}(rh)\mathbf{e}_{k}\right) \right. \\
\displaystyle \left. + E_{h,i}^{-+}(rh-mh)\mathbf{e}_{i}\left(\mathbf{e}_{j}^{-}f_{k}(rh-\mathbf{e}_{j}h)\mathbf{e}_{k}\right)\right]h^{7} = \left\{\begin{array}{cl}
0, &\mbox{if } mh \notin \Omega_{h}\cup\gamma_{h}^{*},\\
f(mh), &\mbox{if } mh \in \Omega_{h}\cup\gamma_{h}^{*},
\end{array} \right.
\end{array}
\end{equation*}
for any discrete function $f$ such that the series converge, and where $E_h^{-+}$ is the discrete fundamental solution to the operator $D_h^{-+}$. Furthermore, $\chi_{\Omega_{h}}$ is the characteristic function of the discrete domain.\par
The discrete octonionic Borel-Pompeiu formula for the exterior domain $\Omega_{h}^{ext}$, which is associated with $\Omega_{h}$, is given by
\begin{equation*}
\begin{array}{c}
\displaystyle \sum\limits_{r\in \mathbb{Z}^{8}}  E_{h}^{-+}(rh-mh)\left( D_h^{+-}f(rh) \right) \chi_{\Omega_{h}^{ext}}(rh) h^{8} \\
\displaystyle +\sum\limits_{rh\in\gamma_{h,2j+1}^{*}}\sum\limits_{j=0}^{7} \left[E_{h}^{-+}(rh-mh)\mathbf{e}_{j}^{+}f(rh+\mathbf{e}_{j}h) + E_{h}^{-+}(rh-mh-\mathbf{e}_{j}h)\mathbf{e}_{j}^{-}f(rh) \right]h^{7} \\
\displaystyle - \sum\limits_{rh\in\gamma_{h,2j}^{*}}\sum\limits_{j=0}^{7} \left[E_{h}^{-+}(rh-mh+\mathbf{e}_{j}h)\mathbf{e}_{j}^{+}f(rh) + E_{h}^{-+}(rh-mh)\mathbf{e}_{j}^{-}f(rh-\mathbf{e}_{j}h) \right]h^{7} \\
\displaystyle - 2\sum_{r\in \mathbb{Z}^{8}} \sum\limits_{s=1}^{7}  \sum_{i\in I_{s}} \sum_{\stackrel{j\in I_{s}}{j\neq i}}^{7} \sum_{\stackrel{k=1}{k\notin I_{s}}}^{7} \left[E_{h,i}^{-+}(rh-mh)\mathbf{e}_{i}\left(\partial_{h}^{+j}\mathbf{e}_{j}^{+}f_{k}(rh)\mathbf{e}_{k}\right) \right. \\
\displaystyle \left. + E_{h,i}^{-+}(rh-mh)\mathbf{e}_{i}\left(\partial_{h}^{-j}\mathbf{e}_{j}^{-}f_{k}(rh)\mathbf{e}_{k}\right)\right]\chi_{\Omega_{h}}^{ext}(rh) h^{8} \\
\displaystyle - 2\sum\limits_{rh\in\gamma_{h,2j+1}^{*}} \sum\limits_{s=1}^{7}  \sum_{i\in I_{s}} \sum_{\stackrel{j\in I_{s}}{j\neq i}}^{7} \sum_{\stackrel{k=1}{k\notin I_{s}}}^{7} \left[E_{h,i}^{-+}(rh-mh)\mathbf{e}_{i}\left(\mathbf{e}_{j}^{+}f_{k}(rh+\mathbf{e}_{j}h)\mathbf{e}_{k}\right) \right. \\
\displaystyle \left. + E_{h,i}^{-+}(rh-mh-\mathbf{e}_{j}h)\mathbf{e}_{i}\left(\mathbf{e}_{j}^{-}f_{k}(rh)\mathbf{e}_{k}\right)\right]h^{7} \\
\displaystyle + 2\sum\limits_{rh\in\gamma_{h,2j}^{*}} \sum\limits_{s=1}^{7}  \sum_{i\in I_{s}} \sum_{\stackrel{j\in I_{s}}{j\neq i}}^{7} \sum_{\stackrel{k=1}{k\notin I_{s}}}^{7} \left[E_{h,i}^{-+}(rh-mh+\mathbf{e}_{j}h)\mathbf{e}_{i}\left(\mathbf{e}_{j}^{+}f_{k}(rh)\mathbf{e}_{k}\right) \right. \\
\displaystyle \left. + E_{h,i}^{-+}(rh-mh)\mathbf{e}_{i}\left(\mathbf{e}_{j}^{-}f_{k}(rh-\mathbf{e}_{j}h)\mathbf{e}_{k}\right)\right]h^{7} = \left\{\begin{array}{cl}
0, &\mbox{if } mh \notin \Omega_{h}^{ext}\cup\gamma_{h}^{*},\\
f(mh), &\mbox{if } mh \in \Omega_{h}^{ext}\cup\gamma_{h}^{*},
\end{array} \right.
\end{array}
\end{equation*}
for any discrete function $f$ such that the series converge, and where $\chi_{\Omega_{h}^{ext}}$ is the characteristic function of the exterior discrete domain. The index sets $I_{s}$, $s=1,\ldots,7$ in these formulae are defined as follows
\begin{equation*}
\begin{array}{cclcclcclccl}
I_{1} & := & \left\{1,2,4\right\}, & I_{2} & := & \left\{1,3,5\right\}, & I_{3} & := & \left\{1,6,7\right\}, & I_{4} & := & \left\{2,3,6\right\}, \\
I_{5} & := & \left\{2,5,7\right\}, & I_{6} & := & \left\{3,4,7\right\}, & I_{7} & := & \left\{4,5,6\right\}.
\end{array}
\end{equation*}
\end{theorem}
\begin{remark}
It is worth underlining that only boundary layer $\gamma_{h}^{*}$ appears explicitly in the discrete Stokes' and Borel-Pompeiu formulae, but because of shifts in the arguments of functions $rh-\mathbf{e}_{j}h$ and $rh+\mathbf{e}_{j}h$, these formulae require indeed two boundary layers: $\gamma_{h}^{*}$ and $\gamma_{h}^{+}$ for the interior setting, and $\gamma_{h}^{*}$ and $\gamma_{h}^{-}$ for the exterior setting. This is another particularity of the discrete case.
\end{remark}\par
Next, by considering a discrete left monogenic function $f$, we immediately obtain the discrete Cauchy formulae for an arbitrary bounded domain in $h\mathbb{Z}^{8}$:
\begin{theorem}\label{Theorem:Cauchy_formulae}
Let $f$ be a discrete left monogenic function with respect to operator $D^{+-}_h$, and let $E_h^{-+}$ be the discrete fundamental solution to the operator $D_h^{-+}$. Then the (interior and exterior) discrete Cauchy formulae for an arbitrary domain $\Omega_h \subset h\mathbb{Z}^{8}$ are given by
\begin{equation*}
\begin{array}{c}
\displaystyle \sum\limits_{rh\in\gamma_{h,2j}^{*}}\sum\limits_{j=0}^{7} \left[E_{h}^{-+}(rh-mh)\mathbf{e}_{j}^{+}f(rh+\mathbf{e}_{j}h) + E_{h}^{-+}(rh-mh-\mathbf{e}_{j}h)\mathbf{e}_{j}^{-}f(rh) \right]h^{7} \\
\displaystyle - \sum\limits_{rh\in\gamma_{h,2j+1}^{*}}\sum\limits_{j=0}^{7} \left[E_{h}^{-+}(rh-mh+\mathbf{e}_{j}h)\mathbf{e}_{j}^{+}f(rh) + E_{h}^{-+}(rh-mh)\mathbf{e}_{j}^{-}f(rh-\mathbf{e}_{j}h) \right]h^{7} \\
\displaystyle - 2\sum_{r\in \mathbb{Z}^{8}} \sum\limits_{s=1}^{7}  \sum_{i\in I_{s}} \sum_{\stackrel{j\in I_{s}}{j\neq i}}^{7} \sum_{\stackrel{k=1}{k\notin I_{s}}}^{7} \left[E_{h,i}^{-+}(rh-mh)\mathbf{e}_{i}\left(\partial_{h}^{+j}\mathbf{e}_{j}^{+}f_{k}(rh)\mathbf{e}_{k}\right) \right. \\
\displaystyle \left. + E_{h,i}^{-+}(rh-mh)\mathbf{e}_{i}\left(\partial_{h}^{-j}\mathbf{e}_{j}^{-}f_{k}(rh)\mathbf{e}_{k}\right)\right]\chi_{\Omega_{h}}(rh) h^{8} \\
\displaystyle - 2\sum\limits_{rh\in\gamma_{h,2j}^{*}} \sum\limits_{s=1}^{7}  \sum_{i\in I_{s}} \sum_{\stackrel{j\in I_{s}}{j\neq i}}^{7} \sum_{\stackrel{k=1}{k\notin I_{s}}}^{7} \left[E_{h,i}^{-+}(rh-mh)\mathbf{e}_{i}\left(\mathbf{e}_{j}^{+}f_{k}(rh+\mathbf{e}_{j}h)\mathbf{e}_{k}\right) \right. \\
\displaystyle \left. + E_{h,i}^{-+}(rh-mh-\mathbf{e}_{j}h)\mathbf{e}_{i}\left(\mathbf{e}_{j}^{-}f_{k}(rh)\mathbf{e}_{k}\right)\right]h^{7} \\
\displaystyle + 2\sum\limits_{rh\in\gamma_{h,2j+1}^{*}} \sum\limits_{s=1}^{7}  \sum_{i\in I_{s}} \sum_{\stackrel{j\in I_{s}}{j\neq i}}^{7} \sum_{\stackrel{k=1}{k\notin I_{s}}}^{7} \left[E_{h,i}^{-+}(rh-mh+\mathbf{e}_{j}h)\mathbf{e}_{i}\left(\mathbf{e}_{j}^{+}f_{k}(rh)\mathbf{e}_{k}\right) \right. \\
\displaystyle \left. + E_{h,i}^{-+}(rh-mh)\mathbf{e}_{i}\left(\mathbf{e}_{j}^{-}f_{k}(rh-\mathbf{e}_{j}h)\mathbf{e}_{k}\right)\right]h^{7} = \left\{\begin{array}{cl}
0, &\mbox{if } mh \notin \Omega_{h}\cup\gamma_{h}^{*},\\
f(mh), &\mbox{if } mh \in \Omega_{h}\cup\gamma_{h}^{*},
\end{array} \right.
\end{array}
\end{equation*}
and
\begin{equation*}
\begin{array}{c}
\displaystyle \sum\limits_{rh\in\gamma_{h,2j+1}^{*}}\sum\limits_{j=0}^{7} \left[E_{h}^{-+}(rh-mh)\mathbf{e}_{j}^{+}f(rh+\mathbf{e}_{j}h) + E_{h}^{-+}(rh-mh-\mathbf{e}_{j}h)\mathbf{e}_{j}^{-}f(rh) \right]h^{7} \\
\displaystyle - \sum\limits_{rh\in\gamma_{h,2j}^{*}}\sum\limits_{j=0}^{7} \left[E_{h}^{-+}(rh-mh+\mathbf{e}_{j}h)\mathbf{e}_{j}^{+}f(rh) + E_{h}^{-+}(rh-mh)\mathbf{e}_{j}^{-}f(rh-\mathbf{e}_{j}h) \right]h^{7} \\
\displaystyle - 2\sum_{r\in \mathbb{Z}^{8}} \sum\limits_{s=1}^{7}  \sum_{i\in I_{s}} \sum_{\stackrel{j\in I_{s}}{j\neq i}}^{7} \sum_{\stackrel{k=1}{k\notin I_{s}}}^{7} \left[E_{h,i}^{-+}(rh-mh)\mathbf{e}_{i}\left(\partial_{h}^{+j}\mathbf{e}_{j}^{+}f_{k}(rh)\mathbf{e}_{k}\right) \right. \\
\displaystyle \left. + E_{h,i}^{-+}(rh-mh)\mathbf{e}_{i}\left(\partial_{h}^{-j}\mathbf{e}_{j}^{-}f_{k}(rh)\mathbf{e}_{k}\right)\right]\chi_{\Omega_{h}}^{ext}(rh) h^{8} \\
\displaystyle - 2\sum\limits_{rh\in\gamma_{h,2j+1}^{*}} \sum\limits_{s=1}^{7}  \sum_{i\in I_{s}} \sum_{\stackrel{j\in I_{s}}{j\neq i}}^{7} \sum_{\stackrel{k=1}{k\notin I_{s}}}^{7} \left[E_{h,i}^{-+}(rh-mh)\mathbf{e}_{i}\left(\mathbf{e}_{j}^{+}f_{k}(rh+\mathbf{e}_{j}h)\mathbf{e}_{k}\right) \right. \\
\displaystyle \left. + E_{h,i}^{-+}(rh-mh-\mathbf{e}_{j}h)\mathbf{e}_{i}\left(\mathbf{e}_{j}^{-}f_{k}(rh)\mathbf{e}_{k}\right)\right]h^{7} \\
\displaystyle + 2\sum\limits_{rh\in\gamma_{h,2j}^{*}} \sum\limits_{s=1}^{7}  \sum_{i\in I_{s}} \sum_{\stackrel{j\in I_{s}}{j\neq i}}^{7} \sum_{\stackrel{k=1}{k\notin I_{s}}}^{7} \left[E_{h,i}^{-+}(rh-mh+\mathbf{e}_{j}h)\mathbf{e}_{i}\left(\mathbf{e}_{j}^{+}f_{k}(rh)\mathbf{e}_{k}\right) \right. \\
\displaystyle \left. + E_{h,i}^{-+}(rh-mh)\mathbf{e}_{i}\left(\mathbf{e}_{j}^{-}f_{k}(rh-\mathbf{e}_{j}h)\mathbf{e}_{k}\right)\right]h^{7} = \left\{\begin{array}{cl}
0, &\mbox{if } mh \notin \Omega_{h}^{ext}\cup\gamma_{h}^{*},\\
f(mh), &\mbox{if } mh \in \Omega_{h}^{ext}\cup\gamma_{h}^{*},
\end{array} \right.
\end{array}
\end{equation*}
which hold for any discrete function $f$ such that the series converge with index sets given by
\begin{equation*}
\begin{array}{cclcclcclccl}
I_{1} & := & \left\{1,2,4\right\}, & I_{2} & := & \left\{1,3,5\right\}, & I_{3} & := & \left\{1,6,7\right\}, & I_{4} & := & \left\{2,3,6\right\}, \\
I_{5} & := & \left\{2,5,7\right\}, & I_{6} & := & \left\{3,4,7\right\}, & I_{7} & := & \left\{4,5,6\right\}.
\end{array}
\end{equation*}
\end{theorem}\par
By using Theorem~\ref{Theorem:Cauchy_formulae}, we can straightforwardly define the associated discrete octonionic Cauchy transforms:
\begin{definition}
For a discrete $l^p$-function $f$, $1\leq p<+\infty$, defined on the boundary layers $\gamma_{h}^{*}$ and $\gamma_{h}^{+}$ the discrete octonionic interior Cauchy transform is defined by
\begin{equation}\label{Cauchy_transforms_interior}
\begin{array}{c}
\displaystyle \mathcal{C}_{\mathbb{O}}^{+}[f](mh) := \displaystyle \sum\limits_{rh\in\gamma_{h,2j}^{*}}\sum\limits_{j=0}^{7} \left[E_{h}^{-+}(rh-mh)\mathbf{e}_{j}^{+}f(rh+\mathbf{e}_{j}h) + E_{h}^{-+}(rh-mh-\mathbf{e}_{j}h)\mathbf{e}_{j}^{-}f(rh) \right]h^{7} \\
\displaystyle - \sum\limits_{rh\in\gamma_{h,2j+1}^{*}}\sum\limits_{j=0}^{7} \left[E_{h}^{-+}(rh-mh+\mathbf{e}_{j}h)\mathbf{e}_{j}^{+}f(rh) + E_{h}^{-+}(rh-mh)\mathbf{e}_{j}^{-}f(rh-\mathbf{e}_{j}h) \right]h^{7} \\
\displaystyle - 2\sum_{r\in \mathbb{Z}^{8}} \sum\limits_{s=1}^{7}  \sum_{i\in I_{s}} \sum_{\stackrel{j\in I_{s}}{j\neq i}}^{7} \sum_{\stackrel{k=1}{k\notin I_{s}}}^{7} \left[E_{h,i}^{-+}(rh-mh)\mathbf{e}_{i}\left(\partial_{h}^{+j}\mathbf{e}_{j}^{+}f_{k}(rh)\mathbf{e}_{k}\right) \right. \\
\displaystyle \left. + E_{h,i}^{-+}(rh-mh)\mathbf{e}_{i}\left(\partial_{h}^{-j}\mathbf{e}_{j}^{-}f_{k}(rh)\mathbf{e}_{k}\right)\right]\chi_{\Omega_{h}}(rh) h^{8} \\
\displaystyle - 2\sum\limits_{rh\in\gamma_{h,2j}^{*}} \sum\limits_{s=1}^{7}  \sum_{i\in I_{s}} \sum_{\stackrel{j\in I_{s}}{j\neq i}}^{7} \sum_{\stackrel{k=1}{k\notin I_{s}}}^{7} \left[E_{h,i}^{-+}(rh-mh)\mathbf{e}_{i}\left(\mathbf{e}_{j}^{+}f_{k}(rh+\mathbf{e}_{j}h)\mathbf{e}_{k}\right) \right. \\
\displaystyle \left. + E_{h,i}^{-+}(rh-mh-\mathbf{e}_{j}h)\mathbf{e}_{i}\left(\mathbf{e}_{j}^{-}f_{k}(rh)\mathbf{e}_{k}\right)\right]h^{7} \\
\displaystyle + 2\sum\limits_{rh\in\gamma_{h,2j+1}^{*}} \sum\limits_{s=1}^{7}  \sum_{i\in I_{s}} \sum_{\stackrel{j\in I_{s}}{j\neq i}}^{7} \sum_{\stackrel{k=1}{k\notin I_{s}}}^{7} \left[E_{h,i}^{-+}(rh-mh+\mathbf{e}_{j}h)\mathbf{e}_{i}\left(\mathbf{e}_{j}^{+}f_{k}(rh)\mathbf{e}_{k}\right) \right. \\
\displaystyle \left. + E_{h,i}^{-+}(rh-mh)\mathbf{e}_{i}\left(\mathbf{e}_{j}^{-}f_{k}(rh-\mathbf{e}_{j}h)\mathbf{e}_{k}\right)\right]h^{7}.
\end{array}
\end{equation}
Analogously, for a discrete $l^p$-function $f$, $1\leq p<+\infty$, defined on the boundary layers $\gamma_{h}^{*}$ and $\gamma_{h}^{-}$ the discrete octonionic exterior Cauchy transform is defined by
\begin{equation}\label{Cauchy_transforms_exterior}
\begin{array}{c}
\displaystyle \mathcal{C}_{\mathbb{O}}^{-}[f](mh) := \displaystyle \sum\limits_{rh\in\gamma_{h,2j+1}^{*}}\sum\limits_{j=0}^{7} \left[E_{h}^{-+}(rh-mh)\mathbf{e}_{j}^{+}f(rh+\mathbf{e}_{j}h) + E_{h}^{-+}(rh-mh-\mathbf{e}_{j}h)\mathbf{e}_{j}^{-}f(rh) \right]h^{7} \\
\displaystyle - \sum\limits_{rh\in\gamma_{h,2j}^{*}}\sum\limits_{j=0}^{7} \left[E_{h}^{-+}(rh-mh+\mathbf{e}_{j}h)\mathbf{e}_{j}^{+}f(rh) + E_{h}^{-+}(rh-mh)\mathbf{e}_{j}^{-}f(rh-\mathbf{e}_{j}h) \right]h^{7} \\
\displaystyle - 2\sum_{r\in \mathbb{Z}^{8}} \sum\limits_{s=1}^{7}  \sum_{i\in I_{s}} \sum_{\stackrel{j\in I_{s}}{j\neq i}}^{7} \sum_{\stackrel{k=1}{k\notin I_{s}}}^{7} \left[E_{h,i}^{-+}(rh-mh)\mathbf{e}_{i}\left(\partial_{h}^{+j}\mathbf{e}_{j}^{+}f_{k}(rh)\mathbf{e}_{k}\right) \right. \\
\displaystyle \left. + E_{h,i}^{-+}(rh-mh)\mathbf{e}_{i}\left(\partial_{h}^{-j}\mathbf{e}_{j}^{-}f_{k}(rh)\mathbf{e}_{k}\right)\right]\chi_{\Omega_{h}}^{ext}(rh) h^{8} \\
\displaystyle - 2\sum\limits_{rh\in\gamma_{h,2j+1}^{*}} \sum\limits_{s=1}^{7}  \sum_{i\in I_{s}} \sum_{\stackrel{j\in I_{s}}{j\neq i}}^{7} \sum_{\stackrel{k=1}{k\notin I_{s}}}^{7} \left[E_{h,i}^{-+}(rh-mh)\mathbf{e}_{i}\left(\mathbf{e}_{j}^{+}f_{k}(rh+\mathbf{e}_{j}h)\mathbf{e}_{k}\right) \right. \\
\displaystyle \left. + E_{h,i}^{-+}(rh-mh-\mathbf{e}_{j}h)\mathbf{e}_{i}\left(\mathbf{e}_{j}^{-}f_{k}(rh)\mathbf{e}_{k}\right)\right]h^{7} \\
\displaystyle + 2\sum\limits_{rh\in\gamma_{h,2j}^{*}} \sum\limits_{s=1}^{7}  \sum_{i\in I_{s}} \sum_{\stackrel{j\in I_{s}}{j\neq i}}^{7} \sum_{\stackrel{k=1}{k\notin I_{s}}}^{7} \left[E_{h,i}^{-+}(rh-mh+\mathbf{e}_{j}h)\mathbf{e}_{i}\left(\mathbf{e}_{j}^{+}f_{k}(rh)\mathbf{e}_{k}\right) \right. \\
\displaystyle \left. + E_{h,i}^{-+}(rh-mh)\mathbf{e}_{i}\left(\mathbf{e}_{j}^{-}f_{k}(rh-\mathbf{e}_{j}h)\mathbf{e}_{k}\right)\right]h^{7}.
\end{array}
\end{equation}
The index sets $I_{s}$, $s=1,\ldots,7$ in these formulae are defined as follows
\begin{equation*}
\begin{array}{cclcclcclccl}
I_{1} & := & \left\{1,2,4\right\}, & I_{2} & := & \left\{1,3,5\right\}, & I_{3} & := & \left\{1,6,7\right\}, & I_{4} & := & \left\{2,3,6\right\}, \\
I_{5} & := & \left\{2,5,7\right\}, & I_{6} & := & \left\{3,4,7\right\}, & I_{7} & := & \left\{4,5,6\right\}.
\end{array}
\end{equation*}
\end{definition}
Comparing the discrete octonionic Cauchy transform introduced here above with the results from the discrete Clifford analysis for bounded domains presented in \cite{Cerejeiras}, it is evident that both results are structurally similar, as expected, but the non-associativity of octonionic multiplication is reflected in the associator terms. Moreover, because the associator terms for domain do not contain the full discrete Cauchy-Riemann operator, but only several partial derivatives, the discrete octonionic monogenicity of $f$ does not necessary lead to vanishing of these terms.\par
We will list now some important properties of the discrete octonionic Cauchy transforms~(\ref{Cauchy_transforms_interior})-(\ref{Cauchy_transforms_exterior}) without proofs, because proving these properties follows along the same pattern as presented in detail in  \cite{CKKS,Cerejeiras,Krausshar_2}:
\begin{corollary}
Let us consider a discrete bounded domain $\Omega_{h}$ and its associated discrete exterior domain $\Omega_{h}^{ext}$, and the three boundary layers $\gamma_{h}^{+}$, $\gamma_{h}^{*}$, and $\gamma_{h}^{-}$. Then the discrete Cauchy transforms~(\ref{Cauchy_transforms_interior})-(\ref{Cauchy_transforms_exterior}) satisfy the following properties:
\begin{itemize}
\item[(i)] The interior and exterior Cauchy transforms have the following mapping properties:
\begin{equation*}
\begin{array}{ll}
\displaystyle \mathcal{C}_{\mathbb{O}}^{+} \colon l^{p}(\gamma_{h}^{+}\cup \gamma_{h}^{*},\mathbb{O}) \to l^{q}(\Omega_{h},\mathbb{O}), & 1 \leq p,q \leq \infty, \\
\displaystyle \mathcal{C}_{\mathbb{O}}^{-} \colon l^{p}(\gamma_{h}^{-}\cup \gamma_{h}^{*},\mathbb{O}) \to l^{q}(\Omega_{h}^{ext},\mathbb{O}), & 1 \leq p < \infty, \frac{8}{7} < q < \infty. \\
\end{array}
\end{equation*}
\item[(ii)] $D_{h}^{+-} \mathcal{C}_{\mathbb{O}}^{+}[f](mh)=0$, $\forall\, mh\in\Omega_{h}\setminus \gamma_{h}^{+}$.
\item[(iii)] $D_{h}^{+-} \mathcal{C}_{\mathbb{O}}^{-}[f](mh)=0$, $\forall\, mh\in\Omega_{h}^{ext}\setminus \gamma_{h}^{-}$.
\end{itemize}
\end{corollary}\par

\section{Discrete octonionic Hardy spaces for bounded domains}\label{Section:Hardy_bounded}

To round off the construction of the basic ingredients of the discrete octonionic function theory for bounded domain based on the Weyl calculus, it is necessary to discuss discrete Hardy spaces and boundary behaviour of discrete octonionic monogenic functions. The construction of discrete Hardy spaces for bounded domain in $\mathbb{R}^{3}$ and $\mathbb{R}^{n}$ has been presented in \cite{Cerejeiras,Cerejeiras_2}. Further, Hardy spaces have also been defined in the context of discrete octonionic setting for half spaces in \cite{Krausshar_2}. These results have indicated that the non-associativity of the octonionic multiplication does not influence the construction of discrete Hardy spaces, Plemelj projections, and discrete extension operators. The reason is that the construction is based on the calculation of the Fourier symbols of the discrete fundamental solution $E_{h}^{-+}$ on the boundary layers. This machinery has been presented in \cite{CKKS}, where discrete Hardy spaces for half-spaces have been introduced for the first time in the context of discrete Clifford analysis. Later, in \cite{Cerejeiras}, these ideas have been extended to the case of bounded domains in $\mathbb{R}^{n}$. Hence, the results from \cite{Cerejeiras} can be directly adapted to the octonionic case. Therefore, we will only shortly outline important ingredients of all constructions and leave all details to \cite{CKKS,Cerejeiras}.\par
The integral representation of the discrete fundamental solution $E_{h}^{-+}$ is given by
\begin{equation*}
E_{h}^{-+}(mh)=\frac{1}{(2\pi)^{8}}\int\limits_{\xi\in[-\frac{\pi}{h},\frac{\pi}{h}]^{8}} \frac{\tilde{\xi}_{-}}{d^{2}}e^{-i\langle mh,\xi \rangle}d\xi, \quad m\in\mathbb{Z}^{8}
\end{equation*}
where $\tilde{\xi}_{-} = \sum\limits_{j=0}^{7} \left(\mathbf{e}_{j}^{+}\xi_{h}^{-j} + \mathbf{e}_{j}^{-}\xi_{h}^{+j}\right)$ and $\xi=\sum\limits_{j=0}^{7}\left(\mathbf{e}_{j}^{+}+\mathbf{e}_{j}^{-}\right)\xi_{j}$. We need to apply the classical 7-D Fourier transform to this discrete fundamental solution. Since this calculation is not affected by the non-associativity, the results carefully carried out in \cite{Cerejeiras} can be carried over directly to our case. Based on these calculations, the discrete Riesz kernels can be introduced now. It is necessary to remark, that the Riesz kernels for arbitrary bounded domains in $\mathbb{R}^{8}$ are obtained by calculating Fourier symbols of the discrete fundamental solution over the boundary layers $\gamma_{h}^{-}$, $\gamma_{h}^{*}$, and $\gamma_{h}^{+}$. For the purpose of presentation, we will follow the approach from \cite{Cerejeiras} and use a mapping $\boldsymbol\Phi$ between the real line and the respective part of the boundary. Note that this approach is typical for the continuous case, and in the discrete setting it has another advantage, that the mapping is realised via $\frac{\pi}{2}k$-rotations and translations of the uniform lattice, because the discrete boundary is always parallel to the coordinate axes.\par
Based on the discussion presented above, the following convolution kernels can be introduced:
\begin{equation*}
\begin{array}{rcl}
H_{i}^{+}f & := & \displaystyle \boldsymbol{\Phi}^{-1}  \mathcal{F}_h^{(7)} \left[\frac{\widetilde{\underline\xi}_{-,i}}{\underline d}\left(e_{i}^+\frac{h\underline d-\sqrt{4+h^2\underline d^{2}}}{2} +e_{i}^-\frac{2}{h\underline d-\sqrt{4+h^2\underline d^{2}}} \right)\right], \\
H_{i}^{-}f & := & \displaystyle \boldsymbol{\Phi}^{-1}  \mathcal{F}_h^{(7)} \left[\frac{\widetilde{\underline\xi}_{-,i}}{\underline d}\left(e_{i}^+\frac{2}{h\underline d-\sqrt{4+h^2\underline d^{2}}}+e_{i}^-\frac{h\underline d-\sqrt{4-h^2\underline d^{2}}}{2}\right)\right],
\end{array}
\end{equation*}
where $\mathcal{F}_h^{(7)}$ denotes the $7$-dimensional Fourier transform which treats the coefficient $m_i h$, associated with the discretisation in the $e_{i}$ direction, as constant, while the term $\widetilde{\underline\xi}_{-,i}$ omits the basis elements $e_i^+,e_i^-$, i.e. $\widetilde{\underline\xi}_{-}=\sum\limits_{j=0,j\neq i}^{7} e_j^+\xi_{-j}^D+e_j^-\xi_{+j}^D$, and $\boldsymbol{\Phi}$ denotes a mapping of individual boundary parts to the real line.\par
\begin{remark}
In contrast, for the case of half-spaces, the Fourier symbols of the discrete fundamental solution $E_{h}^{-+}$ on the layers $m_{7}=-1$, $m_{7}=0$, and $m_{7}=1$ can be calculated explicitly, see again \cite{CKKS}. These symbols are given by
\begin{equation*}
\begin{array}{rcl}
\displaystyle \mathcal{F}_h^{(7)} E_{h}^{-+}(\underline \xi, 0) & = & \displaystyle \frac{\widetilde{\underline{\xi}}_{-}}{\underline{d}\sqrt{4+h^{2}\underline{d}^{2}}} + \left(\mathbf{e}_{7}^{+}-\mathbf{e}_{7}^{-}\right)\left(\frac{1}{2}-\frac{h\underline{d}}{2\sqrt{4+h^{2}\underline{d}^{2}}}\right), \\
\displaystyle \mathcal{F}_h^{(7)} E_{h}^{-+}(\underline \xi, h) & = & \displaystyle \frac{\widetilde{\underline{\xi}}_{-}}{\underline{d}}\left(\frac{2+h^{2}\underline{d}^{2}}{2\sqrt{4+h^{2}\underline{d}^{2}}}-\frac{h\underline{d}}{2}\right) + \mathbf{e}_{7}^{+}\left(\frac{h\underline{d}}{2\sqrt{4+h^{2}\underline{d}^{2}}}-\frac{1}{2}\right) \\
& & \displaystyle - \mathbf{e}_{7}^{-}\left(-\frac{3h\underline{d}+h^{3}\underline{d}^{3}}{2\sqrt{4+h^{2}\underline{d}^{2}}}+\frac{h^{2}\underline{d}^{2}+1}{2}\right) \\
\displaystyle \mathcal{F}_h^{(7)} E_{h}^{-+}(\underline \xi, -h) & = & \displaystyle \frac{\widetilde{\underline{\xi}}_{-}}{\underline{d}}\left(\frac{2+h^{2}\underline{d}^{2}}{2\sqrt{4+h^{2}\underline{d}^{2}}}-\frac{h\underline{d}}{2}\right) - \mathbf{e}_{7}^{-}\left(\frac{h\underline{d}}{2\sqrt{4+h^{2}\underline{d}^{2}}}-\frac{1}{2}\right) \\
& & \displaystyle + \mathbf{e}_{7}^{+}\left(-\frac{3h\underline{d}+h^{3}\underline{d}^{3}}{2\sqrt{4+h^{2}\underline{d}^{2}}}+\frac{h^{2}\underline{d}^{2}+1}{2}\right), \\
\end{array}
\end{equation*}
where $\mathcal{F}_h^{(7)}$ denotes the 7-dimensional discrete Fourier transform, and $\underline{d}^{2}=\frac{4}{h^{2}}\sum\limits_{j=0}^{6}\sin^{2}\left(\frac{\xi_{j}h}{2}\right)$.
\end{remark}
By the help of the convolution kernels, we can now introduce the following operators:
\begin{equation*}
\begin{array}{rcl}
H_+f (mh) & := & \displaystyle \sum_{i=0}^{7} \left[  \sum_{rh\in\gamma_{i}^+} H^+_i(rh-mh)f(rh)\right] h^{7}, \mbox{ for } mh \in \gamma_h^{+}, \\
\\
H_-f (mh) & := & \displaystyle \sum_{i=0}^{7} \left[ \sum_{rh\in\gamma_{i}^-} H^-_i(rh-mh)f(rh) \right] h^{7}, \mbox{ for } mh\in \gamma_h^{-},
\end{array}
\end{equation*}
which satisfy $(H_+)^2 = (H_-)^2=I$. By the help of these operators, the condition for a discrete function $f$ to be a boundary value of a discrete octonionic monogenic function in $\Omega_{h}$ or $\Omega_{h}^{ext}$ can now be formulated as follows:
\begin{equation*}
\begin{array}{lcl}
f(mh) & = & H_+f(mh), \mbox{ for } \gamma_{h}^{+}, \\
f(mh) & = & H_-f(mh), \mbox{ for } \gamma_{h}^{-}.
\end{array}
\end{equation*}
These conditions finish the construction of the discrete octonionic Hardy spaces for bounded domains in $\mathbb{R}^{8}$:
\begin{definition}\label{Definition:Hardy_spaces}
The space of discrete functions $f\in l^{p}(\gamma_{h}^{+},\mathbb{O})$ whose discrete 7D-Fourier transform fulfills $f= H_+f$ on $\gamma_{h}^{+}$ is called the {\itshape interior discrete octonionic Hardy space} and it is denoted by $h_{p,\gamma_{h}^{+}}^{+}$. Analogously, the space of discrete functions $f\in l^{p}(\gamma_{h}^{-},\mathbb{O})$ whose discrete 7D-Fourier transform fulfills $f=  H_-f$ on $\gamma_{h}^{-}$ is called the {\itshape exterior discrete octonionic Hardy space} and it is denoted by $h_{p,\gamma_{h}^{-}}^{-}$.
\end{definition}\par
We finish this section by introducing the discrete Plemelj or Hardy projections, which are written in terms of the operators $H_{+}$ and $H_{-}$ as follows
\begin{equation*}
P_{+}:=\frac{1}{2}\left(I+ H_{+}\right) \mbox{ and } P_{-}:=\frac{1}{2}\left(I+H_{-}\right).
\end{equation*}
By using these discrete projections, an equivalent characterisation of the discrete Hardy spaces can be provided as
\begin{equation*}
f\in h_{p,\gamma_{h}^{+}}^{+} \Longleftrightarrow P_{+}f=f, \mbox{ and } f\in h_{p,\gamma_{h}^{-}}^{-} \Longleftrightarrow P_{-}f=f.
\end{equation*}\par
Finally, we would like to underline, that other operators from the discrete Clifford analysis setting, such as extension operators and trace operators, can easily be defined, because their construction is not influenced by the non-associativity. Hence, the results from \cite{Cerejeiras} can be directly adapted to the octonionic setting by using the results and definitions from \cite{Krausshar_3}. Therefore, we round off the discussion on the discrete octonionic Hardy spaces for bounded domains here.\par

\section{Summary}

In this paper, we have presented the extension of discrete octonionic analysis to arbitrary bounded domains in $\mathbb{R}^{8}$ in the framework of the Weyl calculus approach. At first we have revisited the half-space setting and presented an explicit proof for the discrete octonionic Stokes' formula for the whole space, where we pointed out that the associator appears now, which has been overseen in the previous works. After that, we presented basic setting for the discrete octonionic function theory for half-spaces $h\mathbb{Z}_{+}^{8}$ and $h\mathbb{Z}_{-}^{8}$.\par
To develop a discrete octonionic function theory for bounded domains in $\mathbb{R}^{8}$, we proved explicitly the discrete octonionic Stokes' formula for a cuboid $K_{h}\subset h\mathbb{Z}^{8}$. After that, by using the approach from the discrete Clifford analysis based on characteristic functions, discrete Stokes' formulae are constructed for arbitrary bounded domains in interior and exterior settings. This approach allows a compact and very general presentation of the discrete interior and exterior Borel-Pompeiu formulae, and Cauchy formulae. Next, discrete octonionic Hardy spaces for bounded domains are constructed. Moreover, we underline the points where the non-associativity of octonionic multiplication plays a crucial role, and where it does not affect the results known from the discrete Clifford analysis setting. Thus, the paper finishes the development of the basic fundaments of discrete octonionic analysis and provides a general framework for further research in this field. In particular, solution of boundary value problems in octonionic setting in practical applications typically requires a discretisation of the problem. Hence, the results presented in this paper provide tools towards addressing boundary value problems for PDEs in octonions.\par

\subsection*{Data Availability Statement}

Data sharing not applicable to this article as no datasets were generated or analysed during the current study.

%\begin{acknowledgements}

%\end{acknowledgements}




\end{document}